\documentclass[onecolumn]{autart}
\usepackage{}    
\usepackage[colorlinks,linkcolor=black,anchorcolor=red,citecolor=red]{hyperref}%

\usepackage{graphicx}           
\usepackage{dcolumn}            
\usepackage{bm}                 
\usepackage{bbm}
 \usepackage{subfiles}
 \usepackage[caption=false,font=footnotesize]{subfig}
\usepackage{epsfig}             
\usepackage{times}              
\usepackage[fleqn]{amsmath}
\usepackage{amssymb}
\usepackage{extarrows}
\usepackage{amsfonts}
\usepackage{mathrsfs}
\usepackage{fancyhdr}
\usepackage{color}
\usepackage{enumerate}
\usepackage{epstopdf}
  \graphicspath{{eps/}}
\usepackage [latin1]{inputenc} %
\usepackage{siunitx}
\usepackage{cite} 
 \newenvironment{pf2}{{\it Proof:\enspace}}{\hfill $\blacksquare$\par}

\usepackage{caption}
\usepackage{diagbox}

\newtheorem{theorem}{Theorem}[section] %
\newtheorem{lemma}{Lemma}[section]

\newtheorem{proposition}[theorem]{Proposition}

\allowdisplaybreaks[4]

\begin{document}

\begin{frontmatter}

\title{Input-to-State Stabilization of a Coupled    ODE-PDE System With Time-Varying Coefficients via Composite Boundary Control}

\author{Yongchun Bi$^{1,2}$},
\author{Jun Zheng$^{2}$}\ead{zhengjun2014@aliyun.com},
\author{Guchuan Zhu$^{3}$},
\author{Jiye Zhang$^{1}$}       
\address{$^{1}$State Key Laboratory of Rail
	Transit Vehicle System, Southwest Jiaotong University,
	Chengdu 611756, Sichuan, China\\
	$^{2}$School of Mathematics, Southwest Jiaotong University,
	Chengdu 611756, Sichuan, China\\
	$^{3}${Department of Electrical Engineering, Polytechnique Montr\'{e}al, P.O. Box 6079, Station Centre-Ville, Montreal, QC, Canada H3T 1J4}
}          

\begin{abstract}

This paper proposes a novel composite boundary feedback control law that ensures input-to-state stability (ISS) for a coupled ODE-parabolic PDE system with time-varying coefficients in both subsystems. In controller design, we circumvent the need to directly solve   coupled time-varying parabolic-hyperbolic kernel equations by employing an analytic pre-defined gain function and a time-varying Volterra kernel  function to design the control law explicitly. In stability analysis, to address the simultaneous challenges of Dirichlet boundary disturbances and time-varying coefficients, we employ the square root of a time-varying positive definite matrix and a superlinear function to construct a nonquadratic Lyapunov function for the ODE and a generalized Lyapunov functional  for the PDE, respectively, in the target system, thereby  establishing the ISS in the
$L^2$-norm of the closed-loop system.  Numerical simulations are
presented to illustrate the effectiveness of the proposed
control scheme.
\end{abstract}
\end{frontmatter}
\section{Introduction}
\label{sec:introduction}
\subsection{Literature Review and Research Motivation}
Coupled ODE-PDE systems arise in a wide range of practical applications, such as thermal-electrochemical models of lithium-ion batteries~\cite{tang2017auto}, wellbore-reservoir drilling models~\cite{leobardo2019}, thermal phase-change phenomena~\cite{KOGA2022}, and thermoacoustic Rijke-tube systems~\cite{deandrade2018}.
Stabilization problems for such systems have received much attention in the past decade, and the    backstepping method has emerged as a powerful tool for control design; see~\cite{kang2017auto,Tang2011,wang2023TAC,dimeglio2018,Antonio2010,Krstic2009SCL,Krstic2019PDEODE,deutscher2018auto,Deutscher2021tac,ren2013scl,diMeglio2020auto,Li2014ECOCV,tang2011franklin,Zhou2012IJC,Anfinsen2018auto,Auriol2018auto} for systems with   constant or spatially varying coefficients.
These works have laid the foundation for the analysis and stabilization of coupled ODE-PDE systems under a  general form; see, e.g.,~\cite{vazquez2026} for surveys on backstepping control of PDEs.

Coupled ODE-PDE systems are often inherently time-varying due to changing operating conditions, environmental effects, or time-dependent interconnections. 
In recent years, the stabilization of coupled ODE-PDE systems with time-varying coefficients---including time-varying-delay systems represented as equivalent ODE-PDE cascades with time-varying transport coefficients---has been investigated in~\cite{LinCai2022,CaiYangLiuZhang2021,Krstic2010TAC,BekiarisLiberisKrstic2012}.  Notably, all existing contributions in the literature are restricted to the case where only one of the two subsystems---either the ODE or the PDE---possesses time-varying coefficients.  In particular, in the equivalent ODE-PDE representations of systems with time-varying input delays, the time-varying coefficient takes a specific form determined by the delay and therefore lacks generality. Consequently, the employed method is difficult to apply directly to systems with general time-varying coefficients. In addition, all existing  works are largely  focused on coupled systems involving ODEs and first-order hyperbolic PDEs, or on systems that can be reduced  to such a structure; see, e.g.,~\cite{LinCai2022}. In these settings, the backstepping design naturally yields a set of first-order hyperbolic kernel equations, which can be effectively solved via the method of characteristics and the successive approximations, thus rendering the stabilization problem tractable.

However, the stabilization of coupled ODE-parabolic PDE systems with time-varying coefficients remains largely unexplored. The main obstacle lies in the fact that, for such settings, the standard backstepping transformation leads to kernel equations that are time-dependent and involve both time derivatives and second-order spatial derivatives---a structure that is no longer amenable to characteristic-based solution techniques. 
In this case, the resulting kernel equations may become  intricately coupled with additional parabolic gain equations, leading to a system of gain-kernel equations for which establishing well-posedness represents   a   challenge in  control design process. As a consequence, the stabilization of coupled ODE-parabolic PDE systems with time-varying coefficients  in both subsystems   remains an open challenge.

In practical applications, external disturbances are unavoidable and may enter the system through different channels in practical ODE-PDE systems due to modeling uncertainties, environmental perturbations, actuation imperfections, measurement inaccuracies, and implementation approximations. Input-to-state stability (ISS) can be adopted to quantify how the system state depends on both the initial data and the external inputs.
Within the ISS framework, the stabilization problem for coupled ODE-parabolic PDE systems without involving time-varying coefficients   has been investigated in~\cite{Zhang2019} for  disturbances only acting on the ODE subsystem, and~\cite{Zhang2021,ZhangYangFeng2025,DlalaBenabdallah2026} for in-domain disturbances and/or different  boundary disturbances, respectively.
In addition, input-to-state stabilization has been studied for coupled ODE-hyperbolic PDE systems without involving time-varying coefficients in~\cite{ZhangHaoQiao2019,CoutinhoAndradeCarvalho2024,ZhangWangLi2020}, as well as for systems with time-varying input delays that can be equivalently represented as ODE-transport PDE cascades with a time-varying transport speed in~\cite{cai2018TAC,CaiMengZhangLiu2018}.

Despite these developments, the input-to-state stabilization problem for coupled ODE-PDE systems with general time-varying coefficients---particularly those  not restricted to specific forms induced by time-varying delays---remains open.   The difficulty is especially acute in the presence of Dirichlet boundary disturbances when Lyapunov methods are employed for ISS analysis, since conventional techniques for handling in-domain or Neumann boundary disturbances are not applicable~\cite{Karafyllis2018book,Mironchenko2020,ZhengZhu2026,Zheng2025}.

Notably, the open-loop ISS of coupled ODE-PDE systems with Dirichlet boundary disturbances has been studied in \cite{Dashkovskiy2023,Karafyllis2018book}  in the absence of time-varying coefficients, and in \cite{Dashkovskiy2020IFAC} for the case of time-varying coefficients appearing only in the ODE plant.
However, addressing input-to-state stabilization becomes considerably more challenging for       coupled systems,
particularly when time-varying coefficients are additionally
present.
Consequently, the design of boundary controllers for achieving ISS of coupled ODE-PDE systems with general time-varying coefficients and various disturbances is of significant theoretical and practical importance.

\subsection{Challenges and Contributions}
This paper investigates the input-to-state stabilization problem for a coupled ODE-parabolic PDE system with time-varying coefficients in both plants and bidirectional boundary coupling, where the ODE dynamics are driven by the PDE boundary trace and the PDE boundary condition at the opposite endpoint is coupled to the ODE state, subject to in-domain, Dirichlet boundary, and Robin boundary disturbances.

Compared with  the  existing works, the problem considered in this paper involves two main challenges.

The first challenge lies in the controller design, more precisely, in establishing the solvability of the coupled gain-kernel equations induced by the standard backstepping method.  In the absence of time-varying coefficients, for coupled ODE-parabolic PDE systems with bidirectional coupling through the PDE boundary traces, the induced coupled gain-kernel equations are independent of time and can be solved sequentially through suitable substitutions when the reaction term  of the PDE is absent~\cite{Tang2011,ZhaoXie2014,Benabdallah2020,Dlala2022}. When the reaction term is present, these coupled  gain-kernel  equations can still be treated by successive approximation~\cite{ZhaoXie2014,ayadi2023}, explicit decoupling~\cite{HeXieZhen2017}, or two-step backstepping constructions~\cite{liu2020scl,Benabdallah2022DCD}.  However, in the presence of time-varying coefficients, even for single parabolic PDEs, the backstepping kernel may depend on time, and the corresponding kernel equations contain the time derivatives of the kernel~\cite{Smyshlyaev2005AUT,bao2022ejc,Meurer2009}. In the coupled ODE-PDE setting, a direct application of standard backstepping method would further lead to coupled time-dependent parabolic-hyperbolic equations for the gain and  the kernel function. In this case, the existing techniques for decoupling the equations or solving them by successive approximations are difficult to apply directly, posing a substantial challenge to the solution of the coupled gain-kernel equations  in the controller design for  coupled ODE-PDE systems with time-varying coefficients.

The second challenge lies in the stability analysis. The main difficulty stems from the interplay among Dirichlet boundary disturbances, time-varying coefficients, and ODE-PDE coupling terms. This difficulty is twofold: (i) Dirichlet boundary disturbances generate boundary terms that cannot be handled by standard coercive quadratic  Lyapunov functionals, and their treatment is essentially different from that of Neumann boundary disturbances; (ii) The ODE-PDE coupling and multiple disturbances produce mixed terms  involving the ODE state, the PDE state, and boundary traces, which are difficult to  handle. Within the Lyapunov stability framework, novel techniques,   including nonquadratic Lyapunov functions (NQLFs), generalized Lyapunov functionals (GLFs),  and refined analytical tools, are required to derive  ISS estimates.

To address the first challenge, we develop a composite backstepping design based on an analytic gain function and a time-varying Volterra kernel. The gain is determined first and then used, through the composite transformation, to construct the composite ODE-state feedback gain. The composite ODE-state feedback gain is subsequently incorporated into a decoupled kernel equation to determine the PDE-related Volterra kernel. This sequential construction leads to a composite state-feedback controller and avoids directly solving the induced coupled parabolic-hyperbolic  equations arising from the application of the standard backstepping method.


To address the second challenge,  an NQLF based on  a superlinear function and  the square root of a time-varying positive definite matrix  is employed to  characterize the ODE state, while a GLF is constructed to characterize the PDE state. Refined inequalities are then developed to estimate the mixed ODE-PDE terms, boundary traces, and disturbances, thereby enabling the ISS analysis of the closed-loop system.

The main contributions of this paper are threefold.
\begin{enumerate}
	\item[(i)]  A generic boundary control problem is investigated within the  ISS  framework for a coupled ODE-parabolic PDE system with general time-varying coefficients in both plants and  multiple types of external disturbances,  ensuring broad applicability across different settings.
	
	\item[(ii)]  A novel stabilization strategy is developed for the   system, based on a composite boundary feedback controller that employs an analytic gain function and a time-varying Volterra kernel function. Through  composite transformations, the induced gain equation is decoupled from the kernel equation with a nonlocal boundary term. This decoupling allows the gain function and the kernel function to be solved sequentially, thereby reducing the complexity caused by the coupled gain-kernel equations in the application of the standard backstepping method.
	
	\item[(iii)]  A new Lyapunov-based analytical technique is developed for   ISS analysis of coupled ODE-parabolic PDE systems with time-varying coefficients and different disturbances via a superlinear function and the square root of a time-varying positive definite matrix,  yielding   an NQLF for the ODE   and a  GLF  for the PDE. The proposed approach substantially extends existing Lyapunov-based ISS results from single PDEs   to the coupled ODE-PDE setting with general time-varying coefficients.
\end{enumerate}

It should be emphasized that the present work differs from the existing literature in several aspects.
Compared with existing   studies on the stabilization for coupled ODE-PDE systems,  which typically restrict time-varying coefficients to either the ODE or PDE subsystem~\cite{LinCai2022,CaiYangLiuZhang2021,Krstic2010TAC,BekiarisLiberisKrstic2012}, this paper accommodates time-varying coefficients simultaneously in both subsystems.
Compared with studies on the input-to-state stabilization    for coupled ODE-PDE systems without time-varying coefficients, as well as for systems with time-varying delays that are equivalently represented as ODE-PDE cascades with specific time-varying transport coefficients~\cite{Zhang2019,ZhangWangLi2020,Zhang2021,ZhangYangFeng2025,DlalaBenabdallah2026,ZhangHaoQiao2019,CoutinhoAndradeCarvalho2024,cai2018TAC,CaiMengZhangLiu2018}, this paper addresses input-to-state stabilization for a coupled ODE-PDE system with   general time-varying coefficients and Dirichlet boundary disturbances within the Lyapunov-based framework.
In addition, compared with  the works of \cite{ZhangHaoQiao2019,CoutinhoAndradeCarvalho2024,Zhang2019,ZhangWangLi2020,Zhang2021,ZhangYangFeng2025,DlalaBenabdallah2026,Dashkovskiy2023,Krstic2010TAC,cai2018TAC},  which  rely on standard quadratic Lyapunov functions for the ODEs in the   stability analysis   of coupled ODE-PDE systems, this paper employs      an NQLF for the ODE.
\subsection{Organization and Notation}
In the rest of the paper, we  introduce first some basic notations. In Section~\ref{sec:Problem formulation}, we present the problem formulation. In Section~\ref{controller-design}, a composite state feedback controller is designed by using an analytic gain and a time-varying kernel. The uniform boundedness with respect to (w.r.t.) the time variable and invertibility of backstepping transformations, together with the boundedness of their inverses, are also established. In Section~\ref{stability-analysis}, the closed-loop stability is analyzed by employing
the properties of matrix square roots and constructing a nonquadratic function and a GLF. The main result of this paper, namely the ISS in the $L^2$-norm for the coupled system, is then established.
Numerical {simulations are conducted} in Section~\ref{sec:numerical results} to illustrate the effectiveness of the proposed control scheme. Some conclusions are given in Section~\ref{conclusion}. The proofs of the existence of the kernel solution and the bounded invertibility of the backstepping transformations are provided in the Appendix.

\textit{Notation:} 	Let ${\mathbb{N}_0}:=\{0,1,2,...\}$, ${\mathbb{N}}:=\mathbb{N}_0\setminus\{0\}$, $\mathbb{R}:=(-\infty,+\infty)$, $\mathbb{R}_{\geq 0}:=[0,+\infty)$, $\mathbb{R}_{> 0}:=(0,+\infty)$, and $\mathbb{R}_{\leq 0}:=(-\infty,0]$.

For $m,n,N\in\mathbb N$, $\mathbb R^{m\times n}$ denotes the linear space of $m\times n$ real matrices  equipped with the induced   2-norm \(\|\cdot\|\), and, for notational simplicity, let
$\mathbb R^N:=\mathbb R^{N\times 1}$, which is an $N$-dimensional Euclidean space with the Euclidean norm $|\cdot|$.   
Let $I_N\in\mathbb R^{N\times N}$ denote the identity matrix. For a symmetric matrix $P \in\mathbb R^{N\times N}$, $\lambda_{\min}(P )$ and $\lambda_{\max}(P)$ denote its smallest and largest eigenvalues, respectively. The superscript $\top$ represents matrix transposition. For $E\in\mathbb R^{N\times N}$ and $x\in\mathbb R$, $e^{Ex}$ denotes the matrix exponential defined by
$
e^{Ex}:=\sum_{k=0}^{\infty}\frac{(Ex)^k}{k!}.
$

Let $\Omega$ be a domain (either open or closed) in $\mathbb R^1$ or $\mathbb R^2$. For  a normed linear space $Y$, let $C(\Omega;Y):=\{v:\Omega\to Y\mid v\text{ is continuous on }\Omega\}$. For $i\in\mathbb N\cup\{\infty\}$, let $C^i(\Omega;Y):=\{v:\Omega\to Y\mid v\text{ has continuous derivatives up to order }i\text{ on }\Omega\}$, where $i=\infty$ means that derivatives of all orders are continuous. Let $C^{2,1}(\Omega;\mathbb{R}):=\{v:\Omega \rightarrow\mathbb{R}|~v, v_x, v_{xx}, v_t\in C(\Omega;\mathbb{R})\}$.
For   $p\in[1,+\infty)$, let \(L^p(\Omega;Y):=\{v:\Omega\to Y\mid v\text{ is measurable and }\int_{\Omega}\|v(\xi)\|_Y^p\,\mathrm d\xi<\infty\}\), endowed with the norm \(\|v\|_{L^p(\Omega;Y)}:=\left(\int_{\Omega}\|v(\xi)\|_Y^p\,\mathrm d\xi\right)^{1/p}\) for $v\in L^p(\Omega;Y)$. Let \(L^\infty(\Omega;Y):=\{v:\Omega\to Y\mid v\text{ is measurable and } \operatorname*{ess\,sup}_{\xi\in\Omega}$ $\|v(\xi)\|_Y$ $<\infty\}\), endowed with the norm \(\|v\|_{L^\infty(\Omega;Y)}:=\operatorname*{ess\,sup}_{\xi\in\Omega}\|v(\xi)\|_Y\) for $v\in L^\infty(\Omega;Y)$.
Let the product space $\mathcal H:=\mathbb R^{N}\times L^2((0,1);\mathbb R)$ be equipped with the norm $\|(X,u)\|_{\mathcal H}:=|X|+\|u\|_{L^2((0,1);\mathbb{R})}$ for $(X,u)\in \mathcal{H}$.
Let $\mathscr L(L^2((0,1);\mathbb{R})):=\mathscr L(L^2((0,1);\mathbb{R}),L^2((0,1);\mathbb{R}))$ denote the space of all bounded linear operators on $L^2((0,1);\mathbb{R})$, equipped with the operator norm $\|\mathcal{A}\|_{\mathscr L(L^2((0,1);\mathbb{R})}:=\sup\big\{\|\mathcal{A}v\|_{L^2((0,1);\mathbb{R})}\mid v\in L^2((0,1);\mathbb{R}),\ \|v\|_{L^2((0,1);\mathbb{R})}\leq 1\big\}$ for $\mathcal{A} \in \mathscr L(L^2((0,1);\mathbb{R}) $.

Let $\mathcal {K} := \{\gamma:\mathbb{R}_{\geq 0} \rightarrow \mathbb{R}_{\geq 0}|$$\gamma(0)${$=$}{$0$}, $\gamma$ is continuous, strictly increasing\}, $\mathcal {L}$$:=$$\{\gamma:$$\mathbb{R}_{\geq 0}$$\rightarrow$$\mathbb{R}_{\geq 0}|$ $\gamma$ is continuous, strictly decreasing, {$\lim_{s\rightarrow\infty}\gamma(s)$}{$=$}{$0$}\}, $\mathcal {KL} := \{\beta:$$\mathbb{R}_{\geq 0}$$\times$$\mathbb{R}_{\geq 0}$ $\rightarrow$$\mathbb{R}_{\geq 0}|$${\beta}$ is continuous, $\beta(\cdot,t)\in\mathcal {K}$, $\forall t \in \mathbb{R}_{\geq 0}$; $\beta(s,\cdot) \in \mathcal {L},\forall  s \in {\mathbb{R}_{> 0}}\}$.

Let $
\mathcal{D}_0:=\{(x,y)\mid 0\leq y\leq x\leq 1\}
$ and $
\mathcal{D}:=\{(x,y,t)\mid 0\leq y\leq x\leq 1,\ t\in\mathbb{R}_{\geq 0}\}
$.
For $T\in\mathbb{R}_{>0} $, let $Q_T := (0,1) \times  (0,T)$ and $\overline{Q}_T$$:=$$[0,1]$$\times$$[0,T]$.  Let $Q_\infty$$:=(0,1)\times {\mathbb{R}_{>0}}$ and $\overline{Q}_\infty:=[0,1]\times {\mathbb{R}_{\geq 0}}$.

\section{Problem Statement}\label{sec:Problem formulation}
In this paper, we consider the boundary stabilization problem for a coupled ODE-parabolic PDE system with time-varying coefficients subject to different external disturbances:\begin{subequations}\label{original system}
	\begin{align}
		\!\!\!\!	\dot{X}(t)&=A(t) X(t)+B(t) u(0, t)+D(t),\ t\in \mathbb{R}_{>0},\label{ODE}\\
		\!\!\!\!	u_t(x, t)&=u_{x x}(x, t)\!+\!c(t) u(x, t)\!+\!f(x, t),\ (x,t)\in Q_\infty,\!\!\!\label{PDE1} \\
		\!\!\!\!	u_x(0, t)&=l u(0, t)+C(t) X(t)+d_0(t),\ t\in \mathbb{R}_{>0},\label{PDE2}\\
		\!\!\!\!	u(1, t)&=U(t)+d_1(t),\ t\in \mathbb{R}_{>0},\label{PDE3}\\
		\!\!\!\!	X(0)&=X_0,\	u(x,0)=u_0(x),\  x\in(0,1),\label{PDE4}
	\end{align}
\end{subequations}
where  $X:\mathbb{R}_{\geq 0}\rightarrow\mathbb{R}^{N}$ and $u:\overline{Q}_{\infty}\rightarrow  \mathbb{R}$ are the states of the ODE  and  the parabolic PDE, respectively, $A:\mathbb{R}_{\geq 0}\rightarrow \mathbb{R}^{N\times N}$, $B:\mathbb{R}_{\geq 0}\rightarrow \mathbb{R}^{N}$, and $C: \mathbb{R}_{\geq 0}\rightarrow\mathbb{R}^{1\times N}$ are time-varying matrices,   $c: \mathbb{R}_{\geq 0}\rightarrow \mathbb{R}$ is a time-varying coefficient,  $D: \mathbb{R}_{\geq 0}\rightarrow\mathbb{R}^{N}$ represents internal disturbances of the ODE plant, $f: \overline{Q}_{\infty}\rightarrow  \mathbb{R}$ represents in-domain disturbances of the PDE plant while $d_0 :\mathbb{R}_{\geq 0}\rightarrow \mathbb{R}$ and $d_1:\mathbb{R}_{\geq 0}\rightarrow \mathbb{R}$ represent Robin and Dirichlet boundary disturbances, respectively, $X_0\in\mathbb{R}^{N}$ and $u_0:[0,1]\rightarrow\mathbb{R}$ are the initial data, $U:\mathbb{R}_{\geq 0}\rightarrow \mathbb{R}$ is the control input   to be designed, and $l\in \mathbb{R}_{>0}$ is a constant control gain, which is also  to be determined. 

The control objective is to determine the control law $U(t)$ and the constant  $l$ such that  the coupled ODE-PDE system~\eqref{original system} achieves input-to-state stability (ISS)  w.r.t.  external disturbances $D, f, d_0$, and $d_1$.

To achieve the above control objective, we impose the following assumptions. The time-varying coefficients $A$, $B$, $C$, and $c$ are assumed to be analytic w.r.t. $t$. Moreover, there exists a positive constant $\alpha$ such that , for any $ i\in\mathbb N_0$,
\begin{subequations}
	\begin{align}
		&\sup_{t\in\mathbb{R}_{\geq 0}}\!
		\left\|\frac{\mathrm d^i}{\mathrm dt^i} A(t)\right\|
		\le \alpha^{i+1}i!, \label{eq:analytic-factorial-assumption}\\
		&\sup_{t\in\mathbb{R}_{\geq 0}}
		\left|\frac{\mathrm d^i}{\mathrm dt^i} B(t)\right|+\!	\sup_{t\in\mathbb{R}_{\geq 0}}
		\left|\frac{\mathrm d^i}{\mathrm dt^i} C(t)\right|  +\!	\sup_{t\in\mathbb{R}_{\geq 0}} \left|\frac{\mathrm d^i}{\mathrm dt^i} c(t)\right| \!	 \le \alpha^{i+1}i!.\!\!\label{eq:analytic-factorial-assumption-3}
	\end{align}
\end{subequations}
The disturbances and initial  data satisfy
	$		D \in C(\mathbb{R}_{\ge 0};\mathbb{R}^{N})$,     $d_0,d_1 \in  C^1(\mathbb{R}_{\ge 0};\mathbb{R}), $ $
	f \in C^1(\overline{Q}_{\infty};\mathbb{R})$,  and $u_0 \in C^2((0,1);\mathbb{R})\cap C([0,1];\mathbb{R})$.
In addition, the initial data $u_0$, $X_0$ and boundary disturbances $d_0,d_1$ are assumed to satisfy  certain compatibility conditions, which will be specified later.

Furthermore, we   assume that  there exists an  analytic    time-varying feedback gain $K:\mathbb{R}_{\geq 0}\rightarrow\mathbb{R}^{1\times N}$ satisfying \begin{align}\label{eq:analytic-factorial-assumption-2}
	\sup_{t\in\mathbb{R}_{\geq 0}}
	\left|\frac{\mathrm d^i}{\mathrm dt^i} K(t)\right|
	\le \alpha^{i+1}i!,\  \forall i\in\mathbb N_0,
\end{align}   such that
\begin{equation}
	\dot X(t)=F(t)X(t),\  t\in\mathbb{R}_{\geq0},
	\label{eq:F-system}
\end{equation}
with  \(F(t):=A(t)+B(t)K(t)\)   is uniformly exponentially stable.

\begin{rem}\label{uniform-stability}
	For the linear   system~\eqref{eq:F-system},
	the uniform exponential stability of the zero solution can be determined by any of the following sufficient conditions:
	\begin{itemize}
		\item[(i)] Let \(\Phi(t,s)\) be the state transition matrix of
		\eqref{eq:F-system}. There exist constants \(M\geq1\) and
		\(\alpha>0\) such that
		$
		\|\Phi(t,s)\|
		\leq
		M e^{-\alpha(t-s)}$  for any $ t\geq s\geq0;
		$
		see,
		\cite[Theorem~4.11]{khalil2002nonlinear}.

		\item[(ii)] 
		There exist constants \(\lambda>0\) and
		\(\gamma>0\) such that
		$
		\int_s^t \lambda_{\max}(F(\sigma)+F^\top(\sigma))\,\text{d}\sigma
		\leq
		-\lambda(t-s)+\gamma
		$
		for any $t\geq s\geq0$; 
		see,
		\cite[8.4~Corollary]{Rugh:1996}.
		
		\item[(iii)] $F(t)$ is continuously differentiable and bounded, and its
		all pointwise eigenvalues  satisfy
		for $t\ge0$, $\operatorname{Re}\lambda_i(F(t))\le -\mu$, $i=1,\ldots,N$ with a constant $\mu>0$. Moreover, $\|\dot F(t)\|\le\beta$ with  a constant $\beta>0$;
		see, \cite[8.7~Theorem]{Rugh:1996}.
		
		\item[(iv)] There exists a constant \(\mu_0<0\) such that
		$
		F^\top(t)+F(t)\leq \mu_0 I_N$ for all $t \geq 0;
		$
		see,~\cite[Lemma~1]{Zhou2021auto}.
	\end{itemize}
\end{rem}

\section{Controller Design}\label{controller-design}
To avoid solving the coupled parabolic-hyperbolic gain-kernel equations with time-varying coefficients arising from the standard backstepping control design, we first prescribe an analytic gain function $\theta$. Based on this choice, a composite gain is constructed and then substituted into a kernel equation to determine the time-varying Volterra kernel $k$, which in turn yields the proposed composite controller. One advantage of this construction is that, although the stability analysis of the transformed target system is more involved, it can   decouple the equation for $\theta$ from the kernel function equation for $k$.

Specifically, to design a composite controller, we first define the composite gain
\begin{align}
	\!\!\!	\Gamma(x,t)
	:=
	\theta(x,t)\!-\!\int_0^1 k(x,y,t)\theta(y,t)\,\mathrm dy,\ (x,t)\in\overline{Q}_\infty,\label{Gamma-def}
\end{align}
where
$\theta$  and  $k$ are to be determined subsequently.  The transformations associated with $\theta$ and $k$ are applied  to transform the original system~\eqref{original system} into a desired target system. Then, we define  the  composite boundary feedback control law as
\begin{align}\label{controller-U}
	U(t):=\Gamma(1,t)X(t)\!+\!\!\int_0^1 \!k(1,y,t)u(y,t)\,\mathrm dy,\ t\in\mathbb{R}_{\geq 0}.
\end{align}
\begin{rem}
	Most existing works on the  backstepping stabilization   for systems~\eqref{original system} with the same coupling structure focus primarily on time-invariant coefficients, for which the coupled gain-kernel equations are more tractable.  When the   reaction term $cu$ is absent, $k$ can be written as a function of the gain $\theta$, reducing the remaining equation to a decoupled ODE~\cite{Tang2011,ZhaoXie2014,Benabdallah2020,Dlala2022}. When the reaction term is present, the  coupled gain-kernel equations  can be transformed into integral equations and treated by successive approximations~\cite{ZhaoXie2014,ayadi2023}.  The integral equations arising in this formulation can be subsequently
	decoupled to obtain explicit solutions for the gain and kernel~\cite{HeXieZhen2017}. Alternatively, a two-step backstepping design retains the reaction term in a transitional system, making the kernel equation of the first transformation homogeneous and allowing analytic solutions to be derived~\cite{liu2020scl,Benabdallah2022DCD}.  However, for the coupled ODE-PDE system with time-varying coefficients in both plants considered here, the standard backstepping method would lead to coupled time-dependent  gain-kernel equations, which are difficult to reduce to a decoupled ODE by using the transformations in~\cite{liu2020scl,Benabdallah2022DCD,HeXieZhen2017} or treat directly by successive approximations   as in~\cite{ZhaoXie2014,ayadi2023}. 
\end{rem}
\begin{rem}
	In the application of the standard backstepping method, the transformation is introduced first, and the gain and kernel are then determined simultaneously from the target system. In the proposed design, the analytic gain \(\theta\) is determined first, a composite transformation is then constructed through~\eqref{Gamma-def}, and \(k\) is finally obtained from a decoupled kernel equation. Accordingly, the control law given by~\eqref{controller-U} consists of a PDE-state feedback term generated by \(k\) and an effective ODE-state feedback term with gain \(\Gamma(1,t)\), and is therefore a composite state-feedback control law given by~\eqref{controller-U} induced by the composite transformation.
\end{rem}
\subsection{The Choice of an Analytic Gain Function}
In this section, we choose an analytic gain function $\theta:\overline{Q}_{\infty}\rightarrow \mathbb{R}^{1\times N}$, which will be used in the  control law~\eqref{controller-U} and subsequent transformation. More precisely, let first $K(t)$   satisfy~\eqref{eq:analytic-factorial-assumption-2} and  ensure that    system~\eqref{eq:F-system} is uniformly exponentially  stable, and then define $\theta$ on $\overline Q_\infty$ as
\begin{equation} \theta(x,t) := \sum_{n=0}^{\infty} \frac{x^{2n}}{(2n)!}\mathcal T^nK(t) + \sum_{n=0}^{\infty} \frac{x^{2n+1}}{(2n+1)!}\mathcal T^nL(t), \label{eq:theta-series} \end{equation}
where the differential operator   \(\mathcal T:C^\infty(\mathbb{R}_{\geq 0};\mathbb R^{1\times N}) \to C^\infty(\mathbb{R}_{\geq 0};\mathbb R^{1\times N})\)   is defined by
$
\mathcal T p(t):=\dot p(t)+p(t)J(t)$ with
$J(t):=F(t)-c(t)I_N$ for $p\in C^\infty(\mathbb{R}_{\geq 0};\mathbb R^{1\times N})$,
and  $L(t):=C(t)+lK(t)$.

Throughout this paper, let
\begin{align}\label{def-a-b}
	a:=(1+l)\alpha\ \ \text{and}\ \ 
	b:=8 \alpha+ 2\alpha^2.
\end{align}
Then, we claim that  $\theta$ defined by~\eqref{eq:theta-series} satisfies matching conditions, as stated in the following theorem.
\begin{theorem}
	The analytical gain $\theta$ defined by~\eqref{eq:theta-series} satisfies
	\begin{align}\label{theta-partial}
		\sup_{(x,t)\in\overline Q_\infty} \left|\partial_t^r\theta(x,t)\right| 
		\le 2a b^r r!e^b, \ \forall r\in\mathbb{N}_0
	\end{align}
	and
	\begin{subequations}\label{theta-eq}
		\begin{align}
			\theta_{xx}(x,t)
			&=\theta_t(x,t) + \theta(x,t)(A(t) + B(t)\theta(0,t)- c(t)I_N),\ (x,t)\in \overline Q_\infty,\label{theta-eq-1}\\
			\theta(0,t)&=K(t),\ t \in \mathbb{R}_{\geq 0},\label{theta-boundary1}\\
			\theta_x(0,t)&=C(t)+l\theta(0,t)\label{theta-boundary2},  t \in \mathbb{R}_{\geq 0},\\
			\theta(x,0)  &=  \theta_0(x),\  x\in [0,1], \label{theta-initial}
		\end{align}
		where
		$
		\theta_0(x):=\sum_{n=0}^{\infty} \frac{x^{2n}}{(2n)!}\mathcal T^nK(0) + \sum_{n=0}^{\infty} \frac{x^{2n+1}}{(2n+1)!}\mathcal T^nL(0).
		$
	\end{subequations}
\end{theorem}
\begin{pf2}
	We first show that the series~\eqref{eq:theta-series} defining $\theta$ is well defined.
	Note that the coefficients
	\(A,B,K,C\), and \(c\) satisfy~\eqref{eq:analytic-factorial-assumption}, \eqref{eq:analytic-factorial-assumption-2}, and~\eqref{eq:analytic-factorial-assumption-3}.
	By the Leibniz's formula,~\eqref{eq:analytic-factorial-assumption-2}, and~\eqref{eq:analytic-factorial-assumption-3},  for any $r\in\mathbb{N}_0$, we have
	$\left\|
	\frac{\mathrm d^r}{\mathrm dt^r}(B(t)K(t))
	\right\|
	\leq
	\sum_{j=0}^r
	\binom{r}{j}
	\alpha^{j+1}j!
	\alpha^{r-j+1}(r-j)!
	\leq
	\alpha^2(2\alpha)^r r!$.
	Therefore, using~\eqref{eq:analytic-factorial-assumption}, \eqref{eq:analytic-factorial-assumption-2}, and~\eqref{eq:analytic-factorial-assumption-3} again, we obtain
	\begin{align}
		\left\|\frac{\mathrm d^r}{\mathrm dt^r} J(t)\right\|
		\leq\! \left\|\frac{\mathrm d^r}{\mathrm dt^r} A(t)\right\|
		\!+\! \left\|\frac{\mathrm d^r}{\mathrm dt^r}(B(t)K(t))\right\| \!+\! \left\|\frac{\mathrm d^r}{\mathrm dt^r}c(t)I_N\right\|
		\leq (2\alpha+\alpha^2)(2\alpha)^r r!,\  \forall r\in\mathbb{N}_0.
		\label{eq:M-mu-gamma-bound}
	\end{align}

	For $  n\in\mathbb{N}_0$, let
	$
	E_n(t):=\mathcal T^nK(t).
	$
	Then, by the definition of $\mathcal{T}$, we have $E_{n+1}(t)=\mathcal T(\mathcal T^nK(t))=\frac{\text{d}}{\text{d}t} (\mathcal T^nK(t))+(\mathcal T^nK(t))J(t)$.
	It follows that
	\begin{align}\label{Pn+1}
		E_{n+1}(t)=\dot E_n(t)+E_n(t)J(t).
	\end{align}
	We prove by induction that, for all $n\in\mathbb{N}_0$,
	\begin{align}\label{eq:induction-bound}
		\!\!\!\!\!\!\sup_{t\in\mathbb{R}_{\geq 0}}
		\left|\frac{\text{d}^r}{\text{d}t^r}E_n(t)\right|
		\!\leq
		\alpha(2\alpha)^r(4 \alpha+ \alpha^2)^n (n+r)!,\ \! \forall r\in\mathbb{N}_0.\!\!
	\end{align}
	Indeed, for $n=0$,
	\eqref{eq:analytic-factorial-assumption-2} gives
	$	\sup_{t\in\mathbb{R}_{\geq 0}}|\frac{\mathrm d^r}{\mathrm dt^r}E_0(t)|=\sup_{t\in\mathbb{R}_{\geq 0}}|	\frac{\mathrm d^r}{\mathrm dt^r}K(t)|
	\leq
	\alpha(2\alpha)^r r!$ for any $ r\in\mathbb{N}_0$,
	which yields \eqref{eq:induction-bound}  directly. Suppose that \eqref{eq:induction-bound} holds for any fixed $n\in\mathbb{N}_0$. Differentiating~\eqref{Pn+1} and applying
	the Leibniz's formula gives
	$\frac{\mathrm d^r}{\mathrm dt^r}E_{n+1}(t)
	=
	\frac{\mathrm d^{r+1}}{\mathrm dt^{r+1}}E_n(t)
	+
	\sum_{j=0}^{r}
	\binom{r}{j}
	\frac{\mathrm d^j}{\mathrm dt^j}E_n(t)
	\frac{\mathrm d^{r-j}}{\mathrm dt^{r-j}}J(t)$.
	Using the induction hypothesis~\eqref{eq:induction-bound},~\eqref{eq:M-mu-gamma-bound}, and $r!\sum_{j=0}^{r}\frac{(n+j)!}{j!}
	=\frac{(n+r+1)!}{n+1}$, we obtain
	$  \left|\frac{\mathrm d^r}{\mathrm dt^r}E_{n+1}(t)\right|  
	\leq
	\alpha (2\alpha)^{r}(4 \alpha+ \alpha^2)^{n+1}(n+r+1)!$ for any $ r\in\mathbb{N}_0$.
	Thus \eqref{eq:induction-bound} holds for $n+1$, and hence for all $n\in\mathbb{N}_0$.
	
	Then, using
	$
	(n+r)!\leq 2^{n+r}n!r!
	$ and~\eqref{eq:induction-bound},
	we obtain
	\begin{align}\label{eq-16}
		\sup_{t\in\mathbb{R}_{\geq 0}}
		\left|
		\frac{\mathrm d^r}{\mathrm dt^r}
		\bigl(\mathcal T^nK(t)\bigr)
		\right|
		=
		\sup_{t\in\mathbb{R}_{\geq 0}}
		\left|\frac{\mathrm d^r}{\mathrm dt^r}E_{n}(t)\right|    
		\leq
		\alpha (2\alpha)^{r} (4 \alpha+ \alpha^2)^n(n+r)!   \leq
		\alpha(8 \alpha+ 2\alpha^2)^{n+r}n!r!.
	\end{align}
	
	Similarly, using an   induction argument together with~\eqref{eq-16},  for all $n\in\mathbb{N}_0$, we obtain
	\begin{align}
		\sup_{t\in\mathbb{R}_{\geq 0}}\!
		\left|
		\frac{\mathrm d^r}{\mathrm dt^r}
		\bigl(\mathcal T^nL(t)\bigr)
		\right|
		\!\leq\!
		(1+l)\alpha(8 \alpha+ 2\alpha^2)^{n+r}\!n!r!, \! \ \forall r\in\mathbb{N}_0.
		\label{eq:L-n-L-final-bound}
	\end{align}

	For  all \(n\in\mathbb N_0\),  using \eqref{eq-16},~\eqref{eq:L-n-L-final-bound}, and~\eqref{def-a-b}, we have
	\begin{subequations}\label{equ.16}
		\begin{align}
			\sup_{t\in\mathbb{R}_{\geq 0}} \left| \frac{\mathrm d^r}{\mathrm dt^r} \bigl(\mathcal T^nK(t)\bigr) \right|   \le ab^{n+r}n!r!, \ \forall  r\in\mathbb N_0,\label{16}\\
			\sup_{t\in\mathbb{R}_{\geq 0}}\left| \frac{\mathrm d^r}{\mathrm dt^r} \bigl(\mathcal T^nL(t)\bigr) \right|  \le ab^{n+r}n!r!, \ \forall    r\in\mathbb N_0.\label{17}
		\end{align}
	\end{subequations}
	Consequently, for \(x\in[0,1]\), the inequalities $(2n+1)!\geq (2n)!\geq (n!)^2$ yield
	$\sum_{n=0}^{\infty}
	\left|
	\frac{x^{2n}}{(2n)!}\mathcal T^nK(t)
	\right|
	\le
	a\sum_{n=0}^{\infty}
	\frac{b^n n!}{(2n)!}\le a\sum_{n=0}^{\infty}
	\frac{b^n }{n!}$ and $
	\sum_{n=0}^{\infty}
	\left|
	\frac{x^{2n+1}}{(2n+1)!}\mathcal T^nL(t)
	\right|
	\le
	a\sum_{n=0}^{\infty}
	\frac{b^n n!}{(2n+1)!}\le a\sum_{n=0}^{\infty}
	\frac{b^n }{n!}$.
	Thus, the series~\eqref{eq:theta-series} converges uniformly on
	$\overline Q_\infty$. 	
	
	The same argument applies to the differentiated
	series required for $\theta_t,\theta_x$, and $\theta_{xx}$. Moreover, by~\eqref{eq:theta-series} and~\eqref{equ.16},  we also have
	$	\sup_{(x,t)\in\overline Q_\infty} \left|\partial_t^r\theta(x,t)\right|  \le ab^rr! \sum_{n=0}^{\infty}b^nn! ( \frac{1}{(2n)!} + \frac{1}{(2n+1)!} ) \le 2a b^r r!e^b$,
	which implies~\eqref{theta-partial}.
	
	We now verify that $\theta$ satisfies the required equation~\eqref{theta-eq}. First, the conditions~\eqref{theta-boundary1} and~\eqref{theta-initial} follow directly from the definition of $\theta$. Indeed, evaluating~\eqref{eq:theta-series} at $t=0$ and $x=0$, respectively, yields
	$\theta(x,0)   =  \sum_{n=0}^{\infty} \frac{x^{2n}}{(2n)!}\mathcal T^nK(0)+ \sum_{n=0}^{\infty} \frac{x^{2n+1}}{(2n+1)!}\mathcal T^nL(0)$ and
	$\theta(0,t)  = \mathcal T^0K(t) = K(t)$.
	Thus, \eqref{theta-boundary1} and~\eqref{theta-initial}   are satisfied.
	
	Next, we verify the boundary condition~\eqref{theta-boundary2}.  Indeed, differentiating~\eqref{eq:theta-series} w.r.t. $x$ yields
	$	 \theta_x(x,t)  = \sum_{n=1}^{\infty} \frac{x^{2n-1}}{(2n-1)!}\mathcal T^nK(t) + \sum_{n=0}^{\infty} \frac{x^{2n}}{(2n)!}\mathcal T^nL(t)$. 
	Setting $x=0$ gives
	$
	\theta_x(0,t) = \mathcal T^0L(t) = L(t),
	$
	which, along with the definition of $L(t)$, shows that~\eqref{theta-boundary2} is satisfied.
	
	It remains to verify the equation~\eqref{theta-eq-1}. Differentiating~\eqref{eq:theta-series} twice w.r.t. $x$ and once w.r.t. \(t\), respectively, we have 
	\begin{align}
		\!\!\!\theta_{xx}(x,t)
		&\!=\!\! \sum_{n=0}^{\infty}\! \frac{x^{2n}}{(2n)!}\mathcal T^{n\!+1}K(t) \!+\!\! \sum_{n=0}^{\infty} \!\frac{x^{2n+1}}{(2n\!+\!1)!}\mathcal T^{n\!+1}L(t),\!\label{eq:theta-xx}\\
		\!\!\!\!\theta_t(x,t)
		&\!= \!\!\sum_{n=0}^{\infty} \!\frac{x^{2n}}{(2n)!}
		\frac{\text{d}}{\text{d}t}\bigl(\mathcal T^nK(t)\bigr)\!+ \!\sum_{n=0}^{\infty}\! \frac{x^{2n+1}}{(2n\!+1)!}
		\frac{\text{d}}{\text{d}t}\bigl(\mathcal T^nL(t)\bigr).\label{theta-t}
	\end{align}
	Combining~\eqref{eq:theta-series} and~\eqref{theta-t}, we obtain
	\begin{align*}
	\theta_t(x,t)+\theta(x,t)J(t)
	&= \sum_{n=0}^{\infty} \frac{x^{2n}}{(2n)!} \left( \frac{\text{d}}{\text{d}t}\big(\mathcal T^nK(t)\big) +  (\mathcal T^nK(t))J(t) \right)  \notag\\
	&\quad+  \sum_{n=0}^{\infty} \frac{x^{2n+1}}{(2n+1)!} \left( \frac{\text{d}}{\text{d}t}\big(\mathcal T^nL(t)\big) + (\mathcal T^nL(t))J(t) \right). \end{align*}
	By the definition of $\mathcal T$, one has
	$\frac{\text{d}}{\text{d}t}(\mathcal T^nK(t)) +( \mathcal T^nK(t))J(t)$$  = \mathcal T^{n+1}K(t)$ and
	$\frac{\text{d}}{\text{d}t}(\mathcal T^nL(t))$$ + (\mathcal T^nL(t))J(t)$$=$$\mathcal T^{n+1}L(t)$.
	Therefore, 
	we obtain
	$	\theta_t(x,t)+\theta(x,t)J(t)
	= \sum_{n=0}^{\infty} \frac{x^{2n}}{(2n)!}\mathcal T^{n+1}K(t) + \sum_{n=0}^{\infty} \frac{x^{2n+1}}{(2n+1)!}\mathcal T^{n+1}L(t)$.
	Comparing it with~\eqref{eq:theta-xx}, we obtain~\eqref{theta-eq-1}.
	Consequently, $\theta$ satisfies~\eqref{theta-eq}. 
	Thus, this completes the proof.
\end{pf2}
\begin{rem}
	If $A(t)\equiv A$, $B(t)\equiv B$, $K(t)\equiv K$, and $C(t)\equiv C$ are constant matrices, and $c(t)\equiv c$ is a constant, then
	$
	\mathcal T^nK=KJ^n,
	\mathcal T^nL=LJ^n.
	$
	Therefore,
	$
	\theta(x)
	=
	\sum_{n=0}^{\infty}
	\frac{x^{2n}}{(2n)!}KJ^n
	+
	\sum_{n=0}^{\infty}
	\frac{x^{2n+1}}{(2n+1)!}LJ^n,
	$
	which is equivalent to
		$	\theta(x)
		=
		\begin{pmatrix}
			K & L
		\end{pmatrix}
		e^{Ex}
		\begin{pmatrix}
			I_N \\
			0
		\end{pmatrix}$
	with $E=\begin{pmatrix}
		0 & J\\
		I_N & 0
	\end{pmatrix}$.
	This reduces to the usual matrix exponential construction of the ODE kernel
	in the constant coefficient backstepping design; see, e.g.,~\cite{Tang2011,Krstic2009TAC}.
\end{rem}
\subsection{Determination of the Kernel Equation With a Nonlocal Term}
With the gain function \(\theta\) determined, we can define the kernel \(k\) needed for defining the control law~\eqref{controller-U} and the subsequent transformation by using the  composite gain \(\Gamma\) in~\eqref{Gamma-def}, which is constructed from \(\theta\) and involves a nonlocal term w.r.t. both $\theta$ and $k$. Specifically, let first
$\mu(t)$ be an analytic function satisfying
\begin{subequations}
	\begin{align}
		\sup_{t\in\mathbb{R}_{\geq 0}}
		\left|\frac{\mathrm d^i}{\mathrm dt^i} \mu(t)\right|
		&\le \alpha^{i+1}i!,\ \forall i\in \mathbb{N}_0,\label{mu-ana}\\
		\underline{\mu} :=\inf_{t\in\mathbb{R}_{\geq 0}}\mu(t)&>0.\label{underlinc}
	\end{align}
\end{subequations}
Next, let $k$ be  defined on \(\mathcal D\) as the solution of  the equation:\begin{subequations}\label{kernel-k}
	\begin{align}
		\!\!\!   	k_{x x}(x, y, t)&=k_{y y}(x, y, t)+(c(t)+\mu(t)) k(x, y, t)  +k_t(x, y, t),\ (x,y,t)\in\mathcal{D},\\
		\!\!\!   	\frac{\text{d}}{\text{d} x}(k(x, x, t))&=-\frac{1}{2}(c(t)+\mu(t)),\ (x,t)\in \overline Q_\infty,\\
		\!\!\!  	k_y(x, 0, t)&=lk(x, 0, t)-\!\theta(x,t)\! B(t) + \!\int_0^1 k(x,y,t)\theta(y,t)B(t)\,\mathrm dy =lk(x, 0, t)-\!\Gamma(x, t) B(t),\ (x,t)\!\in\! \overline Q_\infty,\!\! \label{k-3}\\
		\!\!\!   	k(0,0,t)&=0,\  t\in\mathbb{R}_{\geq 0},
	\end{align}
\end{subequations}
where $	\frac{\text{d}}{\text{d} x}(k(x, x, t)):=k_x(x,y,t)|_{y=x}+k_y(x,y,t)|_{y=x}$.

In the  sequel, we   establish the existence and uniqueness of solutions to  equation~\eqref{kernel-k}, as stated in the following theorem, whose proof is provided in Appendix~\ref{kernel-exitence-proof}.
\begin{theorem}\label{kernel-solution-exitence}
	The equation~\eqref{kernel-k} admits a unique solution $k\in C^{\infty}(\mathbb{R}_{\geq 0};C^2(\mathcal{D}_0))$ satisfying   
	\begin{align*}
		\sup_{t\in\mathbb R_{\ge0}}
	\left|\partial_t^i k(x,y,t)\right| \le
	2x\gamma^{2i+3}(i+1)!e^{\gamma^2(x^2-y^2)}
	\left(1+\gamma^2(x^2-y^2)\right)^{i+1}
	\end{align*}
	for every \(i\in\mathbb N_0\),  where $\gamma:=\max\{\alpha,b,2ae^b,2\}$ 
	. In particular, there exists a constant \(\overline{k}>0\) such that
	\( |k(x,y,t)|\le \overline{k}\) for all \((x,y,t)\in\mathcal{D}\).
\end{theorem}

Compared with standard kernel equations for  parabolic systems with time-varying coefficients, the kernel equation in \eqref{kernel-k} contains an additional Volterra-type nonlocal term $\int_0^1 k(x,y,t)\theta(y,t)\text{d}y$ induced by the composite transformation. This term couples the boundary condition with the unknown kernel itself and prevents a direct application of the usual successive approximation arguments. In this work, we extend the applicability of the successive approximation method to a time-dependent kernel equation with nonlocal coupling.
\begin{rem}
	In standard backstepping designs for parabolic PDE  systems, the kernel boundary condition at $y=0$ usually takes a Dirichlet-, Neumann-, or Robin-type~\cite{Smyshlyaev2005AUT,bao2022ejc,Meurer2009}. In contrast, the boundary condition \eqref{k-3} contains both a local nonhomogeneous coupling term and a Volterra-type nonlocal integral term. Although the kernel equation remains linear in $k$, this coupled and nonlocal boundary structure makes the solvability analysis different from that for the time-varying kernel equations commonly considered in~\cite{Smyshlyaev2005AUT,bao2022ejc,Meurer2009}.
\end{rem}
\subsection{Derivation of the Target System}
To address the difficulty arising from the standard backstepping method and the time dependence of the transformations, we have first constructed the gain \(\theta\) and then substituted it into the kernel equation to determine \(k\). Based on these two sequentially obtained transformation components, we now define the composite transformation and show that it maps the original system into a target system suitable for stability analysis. Moreover, both the composite transformation and its inverse are shown to be uniformly bounded in time. This uniform boundedness is essential for transferring the estimates obtained for the target system back to the original system~\eqref{original system}.

More specifically, define the composite transformation  $\mathcal{G}: (X,u) \mapsto  (X,w)$ on $\mathcal{H}$:
\begin{align}\label{composite-transfor}
	\!\!\!\!      (X(t),w(\cdot,t))\!:=\!(\mathcal G[X,u])(\cdot,t)\!:=\!((\mathcal{K}_2\circ\mathcal{K}_1)[X,u])(\cdot,t)\!\!
\end{align}
with
\begin{subequations}
	\begin{align}
		\!\!\!\!\!	(X(t),w(\cdot,t))&:=(\mathcal{K}_2[X,v])(\cdot,t):= \Big(X(t),v(\cdot,t)-\!\!\int_0^x \!k(\cdot,y,t)v(y,t)\,\mathrm dy\Big),\!\!\label{transfor-2}\\
			\!\!\!\!(X(t),v(\cdot,t))&:= (\mathcal K_1 [X,u])(\cdot,t):=(X(t),u(\cdot,t)-\theta(\cdot,t) X(t)),\label{transfor-1}
		\end{align}
	\end{subequations}
	where $v$ denotes the state of the transitional system   given by~\eqref{v-system} in Appendix~\ref{proof-of-target}.
	
	\begin{rem}
		The composition of the transformations~\eqref{transfor-2} and~\eqref{transfor-1} yields $ w(x,t) =u(x,t)-\int_0^x k(x,y,t)u(y,t)\,\mathrm{d}y -\Gamma(x,t)X(t), $ where \(\Gamma(x,t)\) is defined by~\eqref{Gamma-def}. Although this composite transformation has  a similar form to the standard backstepping transformation~\cite{Tang2011,Krstic2009SCL}, its construction follows a different manner. In the proposed controller, $\Gamma$ is first defined through the transformation for the predefined analytical gain $\theta$ (see~\eqref{Gamma-def}), whereas in standard backstepping design, the gain function is introduced as an unknown and subsequently determined jointly with \(k\) from the coupled gain-kernel equations.
	\end{rem}
	
	The following proposition shows that the composite transformation~\eqref{composite-transfor}   maps the original system~\eqref{original system} into a target system  suitable for the subsequent stability analysis.
	\begin{proposition}\label{target-sytem-prop}
		The original system~\eqref{original system} with the  control law  given by~\eqref{controller-U} can be transformed into  the target system\begin{subequations}\label{target-system}
			\begin{align}
				\!\!\!\!\!\!\!	\dot X(t)&=F(t)X(t)+B(t)w(0,t) +D(t),\ t\in \mathbb{R}_{>0},\label{w-ODE}\\
				\!\!	\!\!\!\!\!	w_t(x,t)&=w_{xx}(x,t)-\mu(t)w(x,t)+f(x,t)-\theta(x,t)D(t) -\int_0^x k(x,y,t)f(y,t)\,\mathrm dy+k(x,0,t)d_0(t)\notag\\
				\!\!\!\!\!\!\!	&\quad+\int_0^x k(x,y,t)\theta(y,t)D(t)\,\mathrm dy, \ (x,t)\in Q_\infty,\label{target-2}\\
				\!\!\!\!\!\!\!		w_x(0,t)&=l w(0,t)+d_0(t),\ t\in(0,T),\label{target-3}\\
				\!\!	\!\!\!\!\!		w(1,t)&=d_1(t),\ t\in \mathbb{R}_{>0},\label{target-4}\\
				\!\!	\!\!\!\!\!		X(0)&=X_0, \
				w(x,0) 	= w_0(x), \ x\in (0,1),\label{target-5}
			\end{align}
		\end{subequations}
		where $w_0(x):=v_0(x)-\int_0^x k(x,y,0)v(y,0)\,\mathrm dy$, and $v_0$ is the initial datum of $v$ (see~\eqref{initial-vale} in Appendix~\ref{proof-of-target}).
	\end{proposition}
	
	The  proof of Proposition~\ref{target-sytem-prop} is provided in  Appendix~\ref{proof-of-target}.

In addition, it can be shown that the composite transformation induced by~\eqref{composite-transfor} is invertible and bounded, with a bounded inverse, as stated in the following proposition, whose proof is provided in Appendix~\ref{proof-of-equivent}.
\begin{proposition}\label{equivalent}
	For any fixed $t\in\mathbb{R}_{>0}$, the composite transformation   $ \mathcal G:(X,u)\mapsto  (X,w)$
	is a bijective linear operator, and both $ \mathcal G$  and
	its inverse ${\mathcal{G}}^{-1}:(X,w)\mapsto  (X,u)$ on $\mathcal{H}$,  given by $  (X,u):=\mathcal G^{-1}[X,w]:=(\mathcal{K}^{-1}_1\circ\mathcal{K}_2^{-1})[X,w]$ are  uniformly bounded w.r.t. $t$. More precisely,
	\begin{subequations}
		\begin{align}
			\|(\mathcal G[X,u])(\cdot,t)\|_{\mathcal{H}}
			\leq
			\mathcal M_0
			\|(X(t),u(\cdot,t))\|_{\mathcal{H}},\label{bound-1}\\
			\|(\mathcal G^{-1}[X,w])(\cdot,t)\|_{\mathcal{H}}
			\leq
			\mathcal M_1
			\|(X(t),w(\cdot,t))\|_{\mathcal{H}},\label{bound-2}
		\end{align}
	\end{subequations}
	where
	$
	\mathcal M_0
	:=
	(1+\overline\theta)\Big(1+\frac{\overline k}{\sqrt 2}\Big),
	$
	$
	\mathcal M_1
	:=
	(1+\overline\theta)(1+e^{\overline k}),
	$ and $\overline\theta:= 2ae^b$ 
	is the  bound of \(\theta\) obtained from~\eqref{theta-partial}.
\end{proposition}


To  ensure  the existence of a (classical) solution to the   original system~\eqref{original system} and the target system~\eqref{target-system},   we assume further that the initial data $u_0,X_0,w_0$ and the boundary disturbances $d_0,d_1$ satisfy  the  compatibility conditions:
$	u_{0x}(0) =lu_0(0)+C(0)X_0+d_0(0)$,
$u_0(1) =\int_0^1k(1,y,0)u_0(y)\,\text{d}y+\Gamma(1,0)X_0+d_1(0),$
$	w_{0x}(0)  = lw_0(0)+d_0(0)$,  and
$	w_0(1) =  d_1(0)$.
Note that  classical PDE theory (see, e.g.,~\cite[Chapter  \uppercase\expandafter{\romannumeral4}]{Ladyzenskaja1968}) ensures that  system~\eqref{target-2}-\eqref{target-5}  admits a unique solution $w\in C^{2,1}(Q_T;\mathbb{R})\cap C\left(\overline{Q}_T;\mathbb{R}\right)$.  Moreover,~\cite[Theorem 2.4.1  p.69]{boyce2021} ensures  that the ODE plant~\eqref{w-ODE} admits a unique solution $X\in C^1(\mathbb{R}_{>0};\mathbb{R}^{N})$, and hence, the target system~\eqref{target-system}  admits a unique solution $(X,w)\in C^1((0,T);\mathbb{R}^{N})\times \left(C^{2,1}(Q_T;\mathbb{R})\cap C\left(\overline{Q}_T;\mathbb{R}\right)\right)$ for any $T\in \mathbb{R}_{>0}$.
Furthermore,   we have the following well-posedness result.
\begin{proposition}
	System~\eqref{original system}  under the control law~\eqref{controller-U} admits a unique solution $(X,u)\in C^1((0,T);\mathbb{R}^{N}) \times( C^{2,1}(Q_T;\mathbb{R})\cap C\left(\overline{Q}_T;\mathbb{R}\right))$ for any $T\in \mathbb{R}_{>0}$.
\end{proposition}
\begin{pf2} It suffices to note that  for any fixed $X\in C^1((0,T);\mathbb{R}^{N}) $  the transformation $u\mapsto  w$ determined by  \eqref{composite-transfor} is invertible in  $    C^{2,1}(Q_T;\mathbb{R})\cap C\left(\overline{Q}_T;\mathbb{R}\right) $.
	The proof is similar to Appendix~\ref{proof-of-equivent} and hence omitted.
\end{pf2}

\section{ISS Assessment of the Closed-Loop System}\label{stability-analysis}
This section is devoted to the input-to-state stability (ISS) assessment of  system~\eqref{original system} under the  composite control law~\eqref{controller-U}.
Since the target system~\eqref{target-system} remains coupled and contains time-varying coefficients and Dirichlet  boundary disturbances, standard quadratic Lyapunov functions are not directly applicable. To address these difficulties,  we employ a superlinear function to construct truncation functions, and then, by exploiting the properties of the square root of a time-varying positive-definite matrix, construct  an NQLF for the ODE and a GLF for the PDE.

\subsection{Auxiliary Tools for NQLF-Based  Analysis}
We first consider the basic property of the positive definite square root of the Lyapunov matrix, which will be used to derive the estimate for the ODE plant.
\begin{lemma} \label{lem:H-square-root}
	Let $P\in C^1(\mathbb{R}_{\geq 0};\mathbb R^{N\times N})$ be bounded, symmetric, and   uniformly positive definite; that is $0<h_1I_N\leq P(t)\leq h_2I_N$ with some constants $h_2\geq h_1>0$.  
	Then, the following properties hold.
	\begin{itemize}
		\item[(i)] There exists a unique symmetric and uniformly positive definite matrix $H\in C^1(\mathbb{R}_{\geq 0};\mathbb R^{N\times N})$, such  that
		$H^2(t)=P(t)$.
	Moreover, 	    $H^{-1}(t)$ is symmetric, namely
	\begin{equation} (H^{-1}(t))^\top=H^{-1}(t). \label{eq37} \end{equation}
	\item[(ii)] $H(t)$ and $H^{-1}(t)$ are  bounded,   satisfying
	\begin{subequations}
		\begin{align}
			\sqrt{h_1} I_N &\leq H(t)\leq\sqrt{h_2} I_N, \ \forall  t\in\mathbb{R}_{\geq 0}, \label{eq:H-bounds}\\
			\frac{1}{\sqrt{h_2}}I_N &\leq H^{-1}(t)\leq \frac{1}{\sqrt{h_1}}I_N,\  \forall  t\in\mathbb{R}_{\geq 0}. \label{eq:H-1-bounds} \end{align}
	\end{subequations}
	\item[(iii)] For \begin{equation}
		Y(t):=H(t)X(t),\  Z(t):=X^\top(t)H(t), \label{YZ} \end{equation} it holds that $Z(t)=Y^\top(t)$ and \begin{equation} \!\!\!\!\! \!\!\!\!\! \!\!\!\!\!\!\sqrt{h_1}|X(t)| \leq |Y(t)|=|Z(t)| \leq \sqrt{h_2}|X(t)|,\   \forall t\in \mathbb{R}_{\geq 0}. \label{eq:Y-Z-X-bound} \end{equation} \end{itemize}
		\end{lemma}
		\begin{pf2}
We prove the three assertions in the order stated.

\emph{(i)} Since $P(t)$ is symmetric and uniformly positive definite, there exists an orthogonal matrix $R(t)$ such that  $
P(t)=R(t)$ $\operatorname{diag}\{\lambda_1(t),\ldots,\lambda_N(t)\}R^\top(t)$, for all $t\in \mathbb{R}_{\geq 0}, $ where $\lambda_i(t)>0$, $i=1,\ldots,N$ denotes the $i$-th eigenvalue of $P(t)$. Then, one has $0<h_1\leq \lambda_i(t)\leq h_2$.
Defining
$H(t):=R(t)\operatorname{diag}\{\sqrt{\lambda_1(t)},\ldots,\sqrt{\lambda_N(t)}\}R^\top(t)$
gives $H^2(t)=P(t)$. 
The uniqueness of $H(t)$ follows from the standard uniqueness property of the positive definite matrix $P(t)$. Moreover, since  $P\in C^1(\mathbb{R}_{\geq 0};\mathbb R^{N\times N})$, we have $H\in C^1(\mathbb{R}_{\geq 0};\mathbb R^{N\times N})$. Finally, $H(t)$ is symmetric and uniformly positive definite, and thus is invertible. Since the inverse of a symmetric nonsingular matrix is symmetric, \eqref{eq37} holds.

\emph{(ii)} The eigenvalues of $H(t)$ are $\sqrt{\lambda_i(t)}$, $i=1,\ldots,N$. Therefore,  the fact that $\sqrt{h_1}\leq \sqrt{\lambda_i(t)}\leq \sqrt{h_2}$ implies \eqref{eq:H-bounds}. Since the eigenvalues of $H^{-1}(t)$ are $1/\sqrt{\lambda_i(t)}$, the same bounds yield \eqref{eq:H-1-bounds}.

\emph{(iii)} By \eqref{YZ} and the symmetry of $H(t)$, we have $Z(t)=X^\top(t)H(t)=(H(t)X(t))^\top=Y^\top(t)$. Hence $|Z(t)|=|Y(t)|$. On the other hand, \eqref{eq:H-bounds} gives $\sqrt{h_1}|X(t)|\leq |H(t)X(t)|\leq \sqrt{h_2}|X(t)|$. Combining this with~\eqref{YZ} yields \eqref{eq:Y-Z-X-bound}. The proof is complete. \end{pf2}
\begin{rem}\label{rem-7}
Note that, for $	H(t)$ satisfying Lemma~\ref{lem:H-square-root} (i), $H(t)$ is defined as the unique symmetric  and uniformly positive definite square root  of $P(t)$, denoted by
	$	H(t):=P^{\frac{1}{2}}(t)$.
	\end{rem}
	
	We next introduce a variant of the time-varying Lyapunov equation for the   system~\eqref{eq:F-system} with $F(t)$ as the dynamics matrix.
	\begin{lemma}\label{lemma-2}
If $F\in C  (\mathbb{R}_{\geq 0};\mathbb R^{N\times N}) $ is bounded and system~\eqref{eq:F-system} is  uniformly exponentially stable, then for any   bounded, symmetric, and uniformly positive definite matrix  \(Q\in C(\mathbb{R}_{\geq 0};\mathbb R^{N\times N})\), there exists a
unique  bounded, symmetric, and
uniformly positive definite matrix  \(H\in C^1(\mathbb{R}_{\geq 0};\mathbb R^{N\times N})\) such that
\begin{align}\label{stability-ODE}
	& \dot{H}^{\top}(t)H(t)+ H(t) \dot{H}^{\top}(t)+F^{\top}(t)  H^{\top}(t)H(t)  +H(t) H^{\top}(t)F(t)=-Q(t),\ t\in\mathbb{R}_{\geq 0}.
\end{align}
\end{lemma}
\begin{pf2}
It is clear that  there exists a
unique   bounded, symmetric, and
uniformly positive definite matrix  \(P\in C^1(\mathbb{R}_{\geq 0};\mathbb R^{N\times N})\) satisfying	
\begin{equation}
	\dot P(t)+F^\top(t)P(t)+P(t)F(t)=-Q(t),
	\  \forall t\in\mathbb R_{\geq0},
	\label{eq:P-equation}
\end{equation}
see, \cite[Theorem~4.12]{khalil2002nonlinear}.
By virtue of  Lemma~\ref{lem:H-square-root}, substituting    $P(t)=H^2(t)$ into \eqref{eq:P-equation}, 
we obtain~\eqref{stability-ODE}.
\end{pf2}

In what follows, let  \(P(t)\)    and \(Q(t)\) satisfy the conditions of Lemma~\ref{lem:H-square-root} and  Lemma~\ref{lemma-2}. Thus, there exist constants
$p_1,p_2,q_1,$ $q_2>0$, independent of \(t\), such that
\begin{equation}
\!\!\!  	p_1I_N\leq P(t)\leq p_2I_N,\!\!\!\!\quad
q_1I_N\leq Q(t)\leq q_2I_N,
\  \forall t\in\mathbb R_{\geq0}.\!\!
\label{Q-bound}
\end{equation}
Then, by virtue of  Lemma~\ref{lem:H-square-root} and Lemma~\ref{lemma-2},  we  choose  the unique   bounded, symmetric, and uniformly positive definite square root matrix  \(H\in C^1(\mathbb{R}_{\geq 0};\mathbb R^{N\times N})\) of $P(t)$ such that   \eqref{stability-ODE} is fulfilled.

With the properties of $H(t)$ established in Lemma~\ref{lem:H-square-root}  and  Lemma~\ref{lemma-2}, we now introduce  truncation functions based on a superlinear function to construct  an NQLF and a GLF for the subsequent stability analysis.
\begin{lemma}\label{lem:g and G}
For an  arbitrary constant   $p\in (1,2)$, define truncation functions
\begin{align}\label{gG}
	g(s ):= \begin{cases}s^p, & s \geq 0 \\
		0, & s <0\end{cases}\ \text{and}\
	G(s ):=\int_0^s  g(\tau) \mathrm{d} \tau.
\end{align}
Then, $g$ and $G$ satisfy\begin{subequations}
	\begin{align}
		g(s ) & \geq 0,\quad  g^{\prime}(s) \geq 0,\quad  G(s) \geq 0,  \ \forall s  \in \mathbb{R},\label{gG1} \\
		g(s ) &=G(s )=0, \ \forall s \in \mathbb{R}_{\leq 0},\label{gG2}\\
		G(s )& =\frac{1}{p+1} g(s )s ,\ \forall s  \in \mathbb{R},\label{gG3}\\
		G(s+\tau)&\leq 2^p \left(G(s)+G(\tau)\right),\ \forall s,\tau  \in \mathbb{R}\label{gG4}.
	\end{align}
\end{subequations}
\end{lemma}
\begin{pf2}
The result follows from direct computations.
\end{pf2}


In what follows,   let $Y$ and $Z$  be defined  by~\eqref{YZ} in Lemma~\ref{lem:H-square-root},  and let $g$ and $G$ be defined  by~\eqref{gG} in Lemma~\ref{lem:g and G}. Moreover, define
$
\zeta(x,t)
:=
f(x,t)-\theta(x,t)D(t)-\int_0^x k(x,y,t)f(y,t)\,\mathrm dy +\int_0^x k(x,y,t)\theta(y,t)D(t)\,\mathrm dy+k(x,0,t)d_0(t),
$
which represents the corresponding terms   in~\eqref{target-2} of the target system~\eqref{target-system}. Based on~\eqref{underlinc}, for any $T\in\mathbb{R}_{>0}$,
we further define
$ M:=\max\left\{\|d_1\|_{L^{\infty}(0,T)},\frac{1}{\underline{\mu}}\|\zeta\|_{L^\infty(0,T)}\right\}.
$

\subsection{Estimate for the ODE Plant via an NQLF}

In this subsection, we  derive  the estimate for  the  ODE plant of the target system~\eqref{target-system}, for which a nonquadratic
Lyapunov function is constructed.  More precisely,  instead of adopting a conventional quadratic Lyapunov function, we employ the NQLF $G(|Y|)+G(|Z|)$ for the stability analysis via  the square root matrix $H(t)$ and the nonlinear function $G(a)$.  

\begin{lemma}\label{ODE-system-propo}
For any $T\in\mathbb{R}_{>0}$, the ODE plant~\eqref{w-ODE} of the target system~\eqref{target-system}  has the following estimate:
\begin{align}\label{ODE-estimate}
	 \frac{\text{d}}{\text{d}t}\left(G(|Y|)+G(|Z|)\right) 
	& \leq -(p+1) \underline{\lambda} G(|Y|)  +(4\varepsilon p|H(t)B(t)|+2\varepsilon p) G(|Y|)+\frac{2 \varepsilon^{-p}}{p+1}   \left|H(t) D(t)\right|^{p+1}\notag\\
	&\quad+ \frac{2\varepsilon^{-p}}{p+1}|H(t) B(t)| M^{p+1}+2 \varepsilon^{-p}|H(t) B(t)|G(w(0, t)-M),\ \forall t \in (0,T),
\end{align}
where 	$\underline \lambda:
=
\frac{q_1}{p_2}$ and $\varepsilon$ is an arbitrary positive  constant to be determined later.
\end{lemma}
\begin{pf2}
By the definition of $G$, we obtain
\begin{align}\label{derivates}
	\frac{\text{d}}{\text{d}t}\left(G(|Y|)+G(|Z|)\right)  =  g(|Y|) \frac{\dot{Y}^{\top}   Y}{|Y|}+g(|Z|) \frac{Z   \dot{Z}^{\top}}{|Z|}.
\end{align}
Then,	differentiating $Y$ and $Z$   and using~\eqref{w-ODE} yield
	$		\dot{Y}(t)
	=\dot{H}(t) X(t)+H(t)F( t)X(t)+H(t) D(t)+H(t) B(t) w(0, t)$
and
	$		\dot{Z}(t)
	=  X^{\top}(t)F^{\top}( t) H(t)+D^{\top}(t)   H(t)+X^{\top}(t) \dot{H}(t) +B^{\top}(t) w(0, t) H(t)$.
Substituting $\dot{Y}$ and 	$\dot{Z}$  into~\eqref{derivates}, we obtain
\begin{align}\label{gYgZ2}
 \!\! \frac{\text{d}}{\text{d}t}\left(G(|Y|)+G(|Z|)\right) 
 	&	=  g(|Y|) \frac{1}{|Y|}\Big(X^ { \top }(t)  \left(F^{\top}(t) H^{\top}(t) H(t)+\dot{H}^{\top}(t) H(t)+H(t) H^{\top}(t)F(t)   +H(t) \dot{H}^{\top}(t)\right) X(t) \notag\\
	\!\!	&\quad   +w(0, t) B^{\top}\!(t) H^{\top}\!(t) H(t) X(t)\!+D^{\top}\!(t) H^{\top}\!(t) H(t) X(t)    + \!X^{\!\top}\!(t) H(t) H^{\!\top}\!(t) w(0, \! t) B(t) \!\notag\\
\!\!	&\quad+\!X^{\!\top}\!\! H(t) H^{\!\top}\!(t) D(t)\!\Big).\! \!
\end{align}

Then, for any $t\in(0,T)$, using~\eqref{stability-ODE} in Lemma~\ref{lemma-2}  and the definitions of $Y$ and $Z$ in~\eqref{gYgZ2}, we obtain
\begin{align}\label{gYgZ}
	  \frac{\text{d}}{\text{d}t}\left(G(|Y|)+G(|Z|)\right) &	= g(|Y|) \frac{1}{|Y|}\Big(-X^{\top}(t) Q(t) X(t)+2 X^{\top}(t) H^2(t) D(t) +2 X^{\top}(t) H^2(t) B(t) w(0, t) \Big)\notag\\
	&=  g(|Y|) \frac{1}{|Y|}\left(-Y^{\top}(t)\left(H^{-1}(t)\right)^{\top} Q(t) H^{-1}(t) Y(t)\right)\notag\\
	&\quad+g(|Y|) \frac{1}{|Y|} 2 Y^{\top}\left(H^{-1}(t)\right)^{\top} H^2(t) B(t) w(0, t) \notag\\
	&\quad +g(|Y|) \frac{1}{|Y|} 2 Y^{\top}(t)\left(H^{-1}(t)\right)^{\top} H^2(t) D(t).
\end{align}

We next need to estimate each term on the right-hand side of~\eqref{gYgZ} separately. Indeed, for any $S\in C(\mathbb{R}_{\geq 0}; \mathbb{R}^{N})$,  using
\(|S(t)|=|H(t)H^{-1}(t)S(t)|\le \|H(t)\|\,|H^{-1}(t)S(t)|\),
\(\|H(t)\|^2\le p_2\) by~\eqref{eq:H-bounds}, and~\eqref{Q-bound}, we obtain
\begin{align}\label{S}
	 S^{\top}(t)(H^{-1}(t))^{\top}Q(t)H^{-1}(t)S(t)  
	&=
	(H^{-1}(t)S(t))^{\top}Q(t)(H^{-1}(t)S(t))  \notag\\
	& \ge q_1|H^{-1}(t)S(t)|^2
	\ge \frac{q_1}{\|H(t)\|^2}|S(t)|^2
	\ge \frac{q_1}{p_2}|S(t)|^2.
\end{align}
By~\eqref{gG3} and~\eqref{S}, 	the first term in~\eqref{gYgZ} is changed into
\begin{align}\label{I1}
  g(|Y|) \frac{1}{|Y|}\left(-Y^{\top}(t)\left(H^{-1}(t)\right)^{\top} Q(t) H^{-1}(t) Y(t)\right) 
  \leq-(p+1)\frac{q_1}{p_2} G(|Y|).
\end{align}

The second term and the third term in~\eqref{gYgZ} can be estimated in the same manner. Indeed, applying the Young's inequality with $\varepsilon>0$,~\eqref{eq37}, and~\eqref{gG3}, we obtain
\begin{align}
	 g(|Y|) \frac{1}{|Y|} 2 Y^{\top}(t)\left(H^{-1}(t)\right)^{\top} H^2(t) B(t) w(0, t)  &	=g(|Y|) \frac{2}{|Y|} Y^{\top}(t) H(t) B(t)(w(0, t)-M)\notag\\
	&\quad+g(|Y|) \frac{2}{|Y|} Y^{\top}(t) H(t) B(t) M \notag\\
	&	\leq 2 \varepsilon^{-p}|H(t)B(t)| G (w(0, t )-M)+4 \varepsilon p|H(t)B(t)| G(|Y|) \notag\\
	&\quad+ \frac{2\varepsilon^{-p}}{p+1}|H(t) B(t)| M^{p+1},\label{e1}\\
	  g(|Y|) \frac{1}{|Y|} 2 Y^{\top}(t)\left(H^{-1}(t)\right)^{\top} H^2(t) D(t)  &=  g(|Y|) \frac{2}{|Y|}Y^{\top}(t) H(t) D(t) \notag\\
	&\leq   2\varepsilon p G(|Y|)+ \frac{2\varepsilon^{-p}}{p+1}   \left|H(t) D(t)\right|^{p+1}.\label{e2}
\end{align}

Substituting~\eqref{I1},~\eqref{e1}, and~\eqref{e2} into~\eqref{gYgZ}, we obtain
the desired estimate~\eqref{ODE-estimate} and complete the proof.
\end{pf2}

\subsection{Estimate for the PDE Plant via a GLF}
We next derive an estimate for the PDE plant~\eqref{target-2}-\eqref{target-5}, which will be combined
with the preceding ODE estimate~\eqref{ODE-estimate} to obtain the ISS estimate of the system~\eqref{original system}.
\begin{lemma}\label{w-system-propo}
For any $T\in\mathbb{R}_{>0}$,	the PDE plant~\eqref{target-2}-\eqref{target-5} of the target system~\eqref{target-system} has the following estimate:
\begin{align}\label{PDE-estimate}
 \frac{\text{d}}{\text{d}t}\Big(\int_0^1 G(w-M) \text{d} x +\int_0^1 G(-w-M) \text{d} x\Big) &\leq-l(p+1) G(w(0, t)-M) + \varepsilon p G(w(0, t)-M)\notag\\
	&\quad+\frac{2\varepsilon^{-p}}{p+1}\left|d_0(t)\right|^{p+1}-\underline{\mu}(p+1) \int_0^1 G(w-M) \,\text{d} x\notag\\
	&\quad	-l(p+1) G(-w(0, t)-M) + \varepsilon p G(-w(0, t)-M)\notag\\
	&\quad -\underline{\mu}(p+1) \int_0^1 G(-w-M) \,\text{d} x, \ \forall t \in (0,T).
\end{align}
\end{lemma}
\begin{pf2}
By the definitions of $g$ and $G$, for any $t\in(0,T)$, we have
\begin{align}\label{derivates-PDE}
	 	\frac{\text{d}}{\text{d}t}\Big(\int_0^1 G(w-M) \text{d} x +\int_0^1 G(-w-M) \text{d} x\Big) 
	 =  \int_0^1 g( w-M) w_t \text{d} x-\int_0^1 g(- w-M) w_t \text{d} x.
\end{align}

	We now estimate each term in the right-hand side of~\eqref{derivates-PDE}. Indeed, by~\eqref{target-2},~\eqref{target-3}, and integration by parts, we have
	\begin{align}\label{wt}
		\int_0^1 g( w-M) w_t \text{d} x:=\mathcal{I}_1+\mathcal{I}_2+\mathcal{I}_3+\mathcal{I}_4+\mathcal{I}_5,
	\end{align}
	where
	$	\mathcal{I}_1(t) :=g(w(1, t)-M) w_x(1, t)$,
	$\mathcal{I}_2(t) :=-g(w(0, t)-M)lw(0,t)$,
	$\mathcal{I}_3(t) :=-g(w(0, t)-M)d_0(t)$,
	$\mathcal{I}_4(t) :=	-\int_0^1 g^{\prime}(w-M) w_x^2 \text{d}x$,  and
	$\mathcal{I}_5(t) :=-\int_0^1 g(w-M)\left(\mu(t) w- \zeta\right) \text{d} x$.
	
	It follows from~\eqref{target-4} that
	\begin{align}\label{w leq M}
		|w(1,t)|\leq M,\ \forall t\in(0,T).
	\end{align}
	For $\mathcal{I}_1(t)$, we deduce from~\eqref{gG2} and~\eqref{w leq M}  that $g(w(1,t)-M)=0$ for all $ t\in (0,T)$. Thus, $\mathcal{I}_1(t)=0$ in $(0,T)$.
	For $\mathcal{I}_2(t)$, invoking~\eqref{gG1} and~\eqref{gG3}, we have	
	$ 	\mathcal{I}_2(t)\!
	=  -lg(w(0, t)\!-\!M) (w(0, t)\!-\!M)-l g(w(0, t)-M)M
	\leq  -l(p+1) G(w(0, t)-M)$.
	For $\mathcal{I}_3(t)$, by applying~\eqref{gG3} and the Young's inequality with the same $\varepsilon$  as~\eqref{e2}, we obtain
	$\mathcal{I}_3(t)
	\leq  \varepsilon p G(w(0, t)-M)+\frac{\varepsilon^{-p}}{p+1}\left|d_0(t)\right|^{p+1}$.
	For $\mathcal{I}_4(t)$, by~\eqref{gG1}, we get $\mathcal{I}_4(t)\leq 0$ in $(0,T)$.
	For $\mathcal{I}_5(t)$, using the definition of $M$,~\eqref{gG3}, and~\eqref{underlinc}, we obtain
	$	\mathcal{I}_5(t)
	\leq  -\underline{\mu}(p+1) \int_0^1 G(w-M) \text{d} x$.
	
	Using $\mathcal{I}_1(t)=0$ together with $\mathcal{I}_4(t)\le 0$, and substituting the estimates on $\mathcal{I}_2$, $\mathcal{I}_3$, and $\mathcal{I}_5$ 
	into~\eqref{wt}, we arrive at
	\begin{align}\label{J2}
		\!\! \int_0^1 g( w-M) w_t \text{d} x&\leq 	-l(p+1) G(w(0, t)-M)  + \varepsilon p G(w(0, t)-M)+\frac{\varepsilon^{-p}}{p+1}\left|d_0(t)\right|^{p+1}\notag\\
		\!\!	&\quad-\underline{\mu}(p+1) \int_0^1 G(w-M) \,\text{d} x.
	\end{align}
	
	For the second term in the right-hand side of~\eqref{derivates-PDE},
	note that equation~\eqref{target-system} is linear. Considering  $-w$
	and proceeding in the same way as in~\eqref{J2}, we obtain for any $t\in (0,T)$,
	\begin{align}\label{J3}
		-\int_0^1 g(- w-M) w_t \text{d} x &\leq 	-l(p+1) G(-w(0, t)-M) + \varepsilon p G(-w(0, t)-M)+\frac{\varepsilon^{-p}}{p+1}\left|d_0(t)\right|^{p+1}\notag\\
		&\quad-\underline{\mu}(p+1) \int_0^1 G(-w-M) \,\text{d} x.
	\end{align}
	
	Substituting~\eqref{J2} and~\eqref{J3} into~\eqref{derivates-PDE}, we obtain~\eqref{PDE-estimate}. This completes the proof.
\end{pf2}

\subsection{Estimate for the Closed-Loop System}
We are now ready to combine the estimates for the ODE plant~\eqref{ODE-estimate} and the PDE plant~\eqref{PDE-estimate} with the
bounded invertibility of the transformations to derive an ISS estimate
for the original closed-loop system~\eqref{original system} under the  composite control law~\eqref{controller-U}, as stated in the following theorem.
\begin{theorem}\label{main-result}
	For any $\varepsilon \in \Big(0, \frac{\underline\lambda}{ 2 +4\sqrt{p_2}\alpha}\Big)$, let  $l \geq \max\{\varepsilon^{-1},\varepsilon^{-2}\} \sqrt{p_2}\alpha  + \varepsilon$. Then, under the proposed control law~\eqref{controller-U}, 	the coupled system~\eqref{original system} is ISS w.r.t. disturbances $D$, $f$, $d_0$, and $d_1$, having the following estimate:
	\begin{align*}
		 |X(T)|+ \left\|u(\cdot,T)\right\|_{L^2((0,1);\mathbb{R})} 
		&\leq   e^{-CT}\mathcal{N}_1 \left(\left|X_0\right|+\left\|u_0\right\|_{L^2((0,1);\mathbb{R})},T\right)+\mathcal{N}_2
		\|f\|_{L^\infty(Q_T;\mathbb{R})}+ \mathcal{N}_3 \|d_1\|_{L^\infty((0,T);\mathbb{R})}\notag\\
		&\quad +\mathcal{N}_4
		\|d_0\|_{L^\infty((0,T);\mathbb{R})} + \mathcal{N}_5
		\sup_{t\in(0,T)}|D(t)|,\  \forall T \in\! \mathbb{R}_{>0},
	\end{align*}
	where
	$	C
	:= \min\{
	\underline{\mu},
	\frac{\underline{\lambda}}{2}
	-\varepsilon
	-2\varepsilon\sqrt{p_2}\,\alpha
	\}>0$,
	$\mathcal{N}_1
	:=2\mathcal{M}_0\mathcal{M}_1C_1$,
	$\mathcal{N}_2
	:=\frac{1+\overline{k}}{\underline{\mu}}
	\mathcal{M}_1C_2,$
	$\mathcal{N}_3
	:=C_2\mathcal{M}_1,$
	$\mathcal{N}_4
	:=\left(
	\frac{\overline{k}}{\underline{\mu}}C_2
	+\frac{2C_1}{C\sqrt{\varepsilon}}
	\right)\mathcal{M}_1$, $
	\mathcal{N}_5
	:=\left(
	\frac{1+\overline{k}}{\underline{\mu}}\,
	\overline{\theta}C_2
	+\frac{2C_1\sqrt{p_2}}{C\sqrt{\varepsilon}}
	\right)\mathcal{M}_1$,
	$
	C_1
	:=\frac{\max\{\sqrt{p_2},1\}}
	{\min\{\sqrt{p_1},1\}}$,  and
	$ C_2
	:=\frac{2C_1}{C\sqrt{\varepsilon}}
	\bigl(\sqrt{p_2}\alpha\bigr)^{\frac{1}{2}}
	+\sqrt{2}C_1$.
\end{theorem}
\begin{pf2}
	By Lemma~ \ref{ODE-system-propo} and Lemma~\ref{w-system-propo}, we obtain
	\begin{align}\label{eq-68}
		& \!\frac{\text{d}}{\text{d}t}\Big(G(|Y|)\!+\!G(|Z|)\!+\!\int_0^1\! G(w\!-\!M) \,\text{d} x \!+\!\int_0^1 \!G(-w\!-\!M) \,\text{d} x\Big) \notag\\
		&\leq (p+1)\Big(-\underline\lambda + 2\varepsilon +4\varepsilon  |H(t)B(t)|\Big) G(|Y|)  + \frac{2\varepsilon^{-p}}{p+1} |H(t) B(t)| M^{p+1} + \frac{2\varepsilon^{-p}}{p+1} \big|H(t) D(t)\big|^{p+1}   \notag\\
		&\quad + \Big( - l(p+1)+2\varepsilon^{-p}  |H(t) B(t)|  + \varepsilon p \Big)  G(w(0,t)-M)   + \frac{2\varepsilon^{-p}}{p+1} |d_0(t)|^{p+1}  -\underline{\mu}(p+1) \int_0^1 G(w-M) \,\text{d}x\notag\\
		&\quad+ \Big( - l(p+1)  + \varepsilon p \Big)  G(-w(0,t)-M) -\underline{\mu}(p+1) \int_0^1 G(-w-M) \,\text{d}x,\ \forall t\in(0,T).
	\end{align}
	
	Since
	$	\varepsilon  \in \Big(0, \frac{\underline\lambda}{ 2 +4\sqrt{p_2}\alpha}\Big)$ and
	$	l \geq \max\{\varepsilon^{-1},\varepsilon^{-2}\} \sqrt{p_2}\alpha  + \varepsilon$, by~\eqref{eq:H-bounds} and~\eqref{eq:analytic-factorial-assumption}, we obtain $- l(p+1)+2\varepsilon^{-p}  |H(t) B(t)|  + \varepsilon p< 0$ and $\underline\lambda - 2\varepsilon -4\varepsilon  |H(t)B(t)|>0$ whenever $p\in(1,2)$.
	
	It follows from~\eqref{eq-68} and the Gronwall's inequality that
	\begin{align}\label{upper-bound}
		&G(|Y(T)|)+G(|Z(T)|)+\int_0^1 G(w(x,T)-M) \,\text{d} x +\int_0^1 G(-w(x,T)-M) \,\text{d} x\notag\\
		&\leq e^{-C(p+1)T}\bigg(G(|Y(0)|)+G(|Z(0)|)+\int_0^1 G(w_0-M) \,\text{d} x  +\int_0^1 G(-w_0-M) \,\text{d} x\bigg)
		\notag\\
		&\quad+\frac{2\varepsilon^{-p}}{p+1}
		\int_0^T e^{-C(p+1)(T-s)} \left(
		|H(s)B(s)|M^{p+\!1}
		\!\!+\!|H(s)D(\!s)|^{p+\!1}
		\!\!+\!|d_0(\!s)|^{p+\!1}
		\right)\!\!\,\mathrm ds,\!
	\end{align}
	where $C:=\min\{\underline{\mu},\frac{\underline\lambda}{2} - \varepsilon -2 \varepsilon \sqrt{p_2}\alpha \}$.
	
	From~\eqref{upper-bound}, we get
	\begin{align*}
		  G(|Y(T)|)+\int_0^1 G(|w(x,T)|-M) \,\text{d} x  
		 &\leq 2e^{-C(p+1)T}(G(|Y(0)|)+\int_0^1 G(|w_0|) \,\text{d} x )\notag\\
		&\quad
		+\frac{2\varepsilon^{-p}}{p+1}
		\int_0^T e^{-C(p+1)(T-s)}
		\left(
		|H(s)B(s)|M^{p+1} +|H(s)D(s)|^{p+1}
		\right.\notag\\
		&\left.\quad+|d_0(s)|^{p+1}
		\right)\,\mathrm ds,
	\end{align*}
	which, along with the definitions of $g$ and $G$, as well as~\eqref{gG3} and~\eqref{gG4}, implies that
	\begin{align}\label{eqY}
	 |Y(T)|^{{ {p}}+1}     + \int_0^1|w(x,T)|^{{ {p}}+1} \text{d} x &\leq e^{-C(p+1)T} 2^{p+1} \left( |Y(0)|^{p+1}  + \int_0^1 |w_0 |^{p+1}   \text{d}x \right)\notag\\
		&\quad+2^{p+1}\varepsilon^{-p}
		\int_0^T e^{-C(p+1)(T-s)}
		\left(
		|H(s)B(s)|M^{p+1}
		\right.\notag\\
		&\left.\quad+|H(s)D(s)|^{p+1}
		+|d_0(s)|^{p+1}
		\right)\,\mathrm ds   +2^p M^{p+1}.
	\end{align}
	
	Taking both sides of~\eqref{eqY} to the power $\frac{1}{p+1}$ and using~\eqref{eq:H-bounds} and~\eqref{eq:Y-Z-X-bound}, we obtain
	\begin{align}\label{X+w}
		 |X(T)|+\|w(\cdot,T)\|_{L^{p+1}((0,1);\mathbb{R})} 
		&\leq 2e^{-CT} C_1 \left( |X(0)|   +  \|w_0 \|_{L^{p+1}((0,1);\mathbb{R})}   \right)+2^{\frac{p}{p+1}}C_1 M\notag\\
		&\quad+2C_1\varepsilon^{-\frac{p}{p+1}}
		\int_0^T e^{-C(T-s)}
		\Big(
		|H(s)B(s)|^{\frac{1}{p+1}}M  +|H(s)D(s)|
		+|d_0(s)|
		\Big)\,\mathrm ds
		\notag\\
		&\leq 2e^{-CT} C_1 \left( |X(0)|   +  \|w_0 \|_{L^{p+1}((0,1);\mathbb{R})}   \right)+2^{\frac{p}{p+1}}C_1M\notag\\
		&\quad+\frac{2C_1}{C}\varepsilon^{-\frac{p}{p+1}}
		\Big(
		(\sqrt{p_2}\alpha)^{\frac{1}{p+1}}M
		+\sqrt{p_2}\sup_{t\in(0,T)}|D(t)|+\|d_0\|_{L^\infty((0,T);\mathbb{R})}
		\Big)
		,
	\end{align}
	where $C_1:=\frac{\max\{\sqrt{p_2},1\}}{\min\{\sqrt{p_1},1\}}$.
	
	According to the definition of $M$, we have
	$M	 \leq \|d_1\|_{L^{\infty}((0,T);\mathbb{R})}+\frac{\overline k}{\underline{\mu}} \|d_0\|_{L^{\infty}((0,T);\mathbb{R})}
	+\frac{1+\overline k}{\underline{\mu}} \big(\|f\|_{L^\infty(Q_T;\mathbb{R})}+\overline \theta \sup_{t\in(0,T)}|D(t)|\big)$.
	Then, substituting it into~\eqref{X+w} and  letting $p\to 1$, we have
	\begin{align*}
		 |X(T)|+ \left\|w(\cdot,T)\right\|_{L^2((0,1);\mathbb{R})} 
		&\leq 2 e^{-CT}C_1\left(\left|X_0\right|+\left\|w_0\right\|_{L^2((0,1);\mathbb{R})},T\right) + \left(
		\frac{1+\overline{k}}{\underline{\mu}}\overline{\theta}C_2
		+\frac{2C_1}{C\sqrt{\varepsilon}}\sqrt{p_2}
		\right)
		\sup_{t\in(0,T)}|D(t)|\notag\\
		&\quad + \frac{1+\overline{k}}{\underline{\mu}}C_2
		\|f\|_{L^\infty(Q_T;\mathbb{R})}+ C_2 \|d_1\|_{L^\infty((0,T);\mathbb{R})}
		\notag\\
		&\quad +\left(\frac{\overline{k}}{\underline{\mu}}C_2
		+\frac{2C_1}{C\sqrt{\varepsilon}} \right)
		\|d_0\|_{L^\infty((0,T);\mathbb{R})} , \ \forall T \in \mathbb{R}_{>0},
	\end{align*}
	where $C_2:=
	\frac{2C_1}{C\sqrt{\varepsilon}}
	(\sqrt{p_2}\alpha)^{\frac{1}{2}}
	+\sqrt{2}C_1$.
	
	By virtue of Proposition~\ref{equivalent},
	we obtain the estimate for the original system~\eqref{original system}.
	The proof of Theorem~\ref{main-result} is complete.
\end{pf2}

\section{Numerical Experiments}\label{sec:numerical results}
In this section, numerical simulations are presented to illustrate the effectiveness of the proposed composite control law~\eqref{controller-U}. The closed-loop responses are shown under multiple disturbances, including the in-domain disturbances, Dirichlet boundary disturbances, and Robin boundary disturbances.

The ODE coefficients are selected as
$A(t) = \begin{pmatrix}
	-1+3\cos^2 t & 1-3\sin t \cos t\\
	-1-3\sin t \cos t & -1+3\sin^2 t
\end{pmatrix}$, 
$
B(t)  = \begin{pmatrix}\cos t&  -\sin t\end{pmatrix}^\top,$
$K(t) = \begin{pmatrix}-2.5\cos t & 2.5\sin t\end{pmatrix},$
$C(t)  = \begin{pmatrix}0.5+0.1\sin t & 0.3+0.1\cos t\end{pmatrix}$,
and the PDE reaction coefficient is \(c(t)=1.2\pi+0.1\sin t\). For the kernel computation, we use \(\mu(t)=0.55+0.05\cos t\).
With \(F(t)=A(t)+B(t)K(t)\), direct calculation yields
$
F^\top(t)+F(t)=-2I_2+B(t)B^\top(t)\leq -I_2,
$
which satisfies the condition presented in Remark~\ref{uniform-stability}(iv), and hence system~\eqref{eq:F-system} without external input is uniformly exponentially stable.

The matrix \(Q(t)\) is selected as
$
Q(t)=
\begin{pmatrix}
	2q_h(t) & 0\\
	0 & 2q_h(t)
\end{pmatrix},
$
where
$
q_h(t)
=  \bigl(0.15+0.01\sin(0.8t)\bigr)
\bigl(0.15+0.01\sin(0.8t)-0.008\cos(0.8t)\bigr).
$
The matrix \(H(t)\) is selected as
$H(t)=
\begin{pmatrix}
	0.15+0.01\sin(0.8t) & 0\\
	0 & 0.15+0.01\sin(0.8t)
\end{pmatrix}$,
which satisfies \(0.14I_2\leq H(t)\leq0.16I_2\).
For \(H(t)\) and \(Q(t)\), numerical evaluation gives
$
\underline{\lambda}
=\frac{q_1}{p_2}
\approx 1.5156>0.
$
Moreover, since \(\sup_{t\in\mathbb{R}_{\geq0}}|H(t)B(t)|=\sqrt{p_2}\alpha=0.16\), 
\(\varepsilon=0.56\in(0,\underline{\lambda}/(2+4\sqrt{p_2}\alpha))\), and the condition for $\varepsilon$ in Theorem~\ref{main-result} is satisfied. For this choice of \(\varepsilon\), the lower bound for \(l\) is
\(\max\{\varepsilon^{-1},\varepsilon^{-2}\}\sqrt{p_2}\alpha+\varepsilon\approx1.0702\). Therefore, \(l=4\) satisfies the required condition for $l$ in Theorem~\ref{main-result}.

The initial condition is
$	u_0(x) =-0.5\phi(x+1)^2(x+4)(x-0.5)(x+0.2),
X_0 =\phi(0.02,-0.02)^\top,$
where \(\phi\in\{0.5,2\}\) represents amplitudes of initial data. 

For disturbances, we take
$	D(t) =\delta\begin{pmatrix}0.03\sin(3t)\\0.04\cos(3t)\end{pmatrix}$,
$f(x,t) =\sigma\sin(2(t+x))$,
$d_0(t) =\sigma_0\sin(2t-\pi)$ and
$ d_1(t) =\sigma_1\sin(2(t-\pi))$
with \(\delta=\sigma\in\{0,2,4\}\) and \(\sigma_0=\sigma_1\in\{0,1,3\}\).

\begin{figure}[!htbp]
	\centering
	\subfloat[]{\includegraphics[scale=0.31]{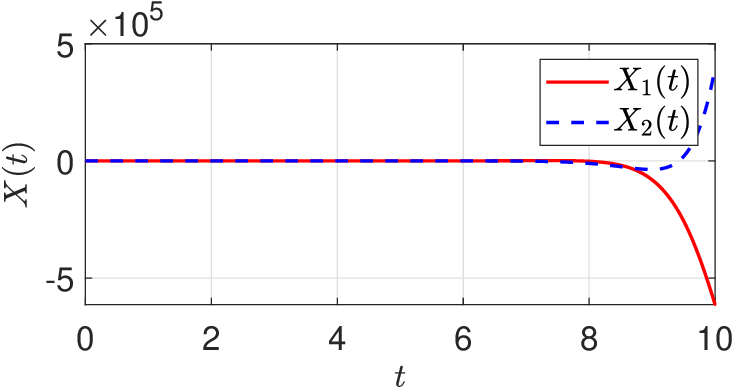}}    \hspace{0.5cm}
	\subfloat[]{\includegraphics[scale=0.31]{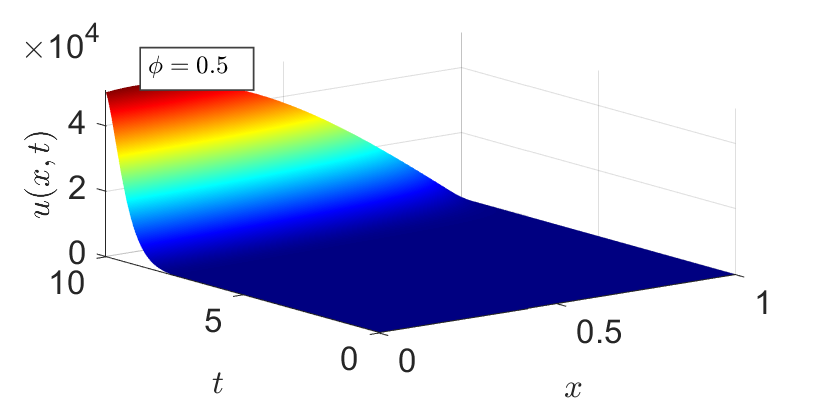}}\\
	\subfloat[]{\includegraphics[scale=0.31]{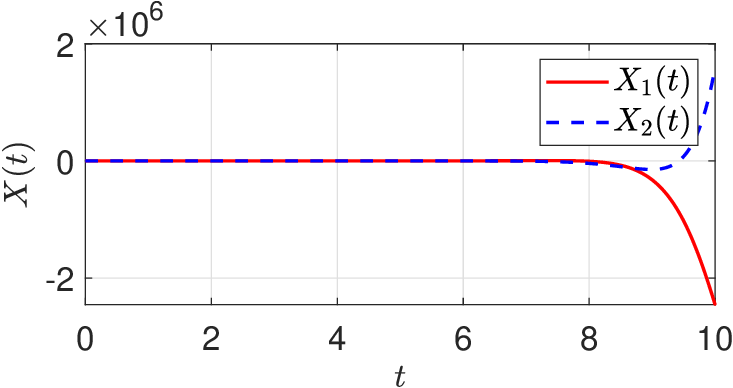}}    \hspace{0.5cm}
	\subfloat[]{\includegraphics[scale=0.31]{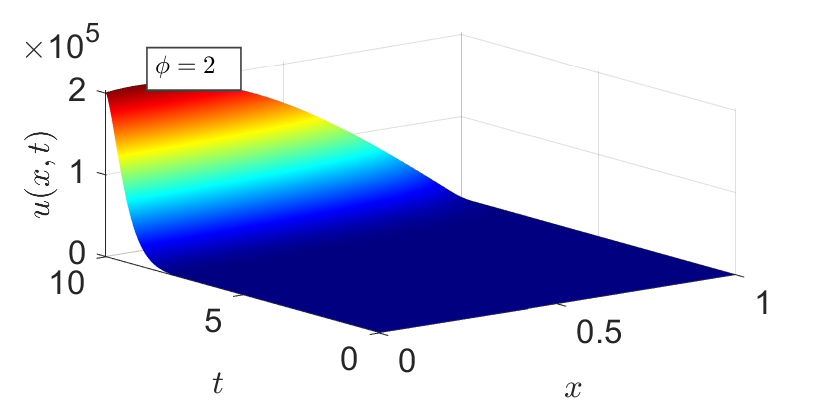}}
	\caption{Evolution of states  in open loop with different initial data: (a)~ODE state when $\phi=0.5$; (b) PDE state when $\phi=0.5$; (c) ODE state when $\phi=2$; (d) PDE state when $\phi=2$. 
	}
	\label{fig:open-loop}
\end{figure}

\begin{figure}[!htbp]
	\centering
	\subfloat[]{\includegraphics[scale=0.31]{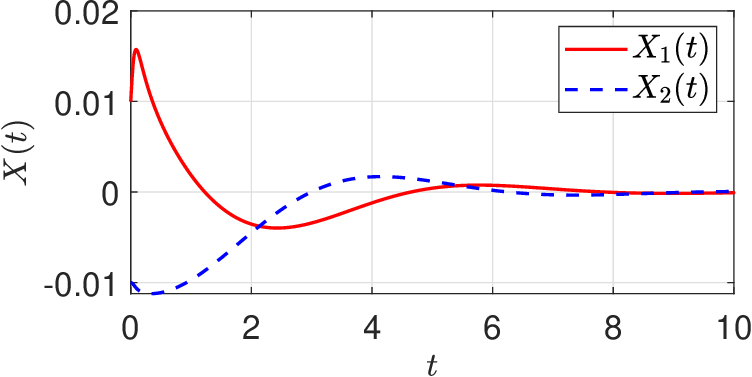}}   \hspace{0.5cm}
	\subfloat[]{\includegraphics[scale=0.31]{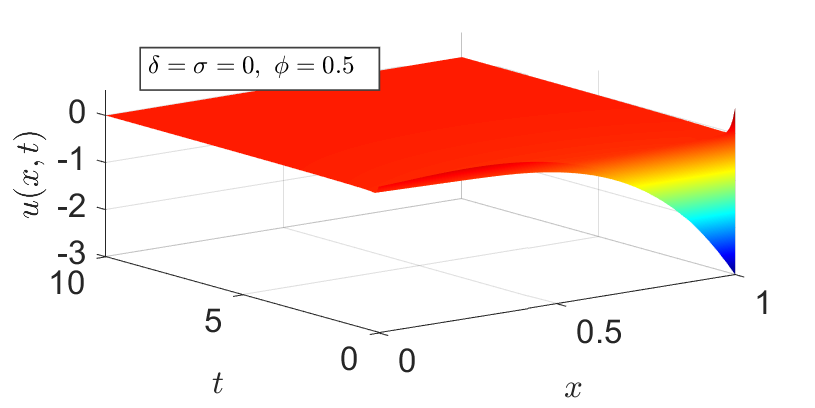}}\\
	\subfloat[]{\includegraphics[scale=0.31]{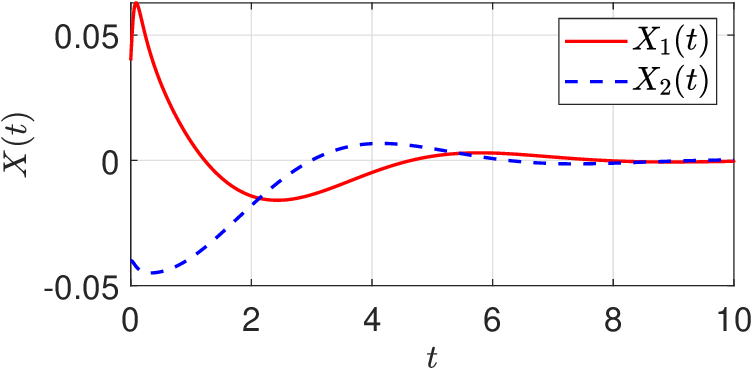}}   \hspace{0.5cm}
	\subfloat[]{\includegraphics[scale=0.31]{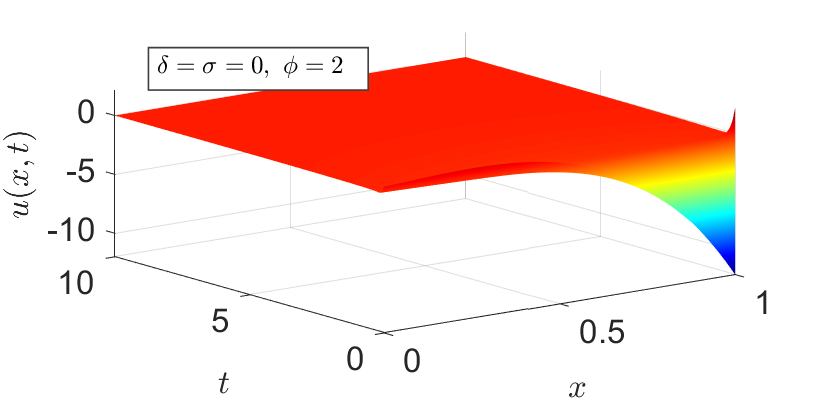}}
	\caption{Evolution of  states   in   closed loop  with different initial data in the absence of disturbances: (a) ODE state when $\phi=0.5$; (b) PDE state when $\phi=0.5$; (c) ODE state when $\phi=2$; (d) PDE state {when} $\phi=2$.}
	\label{fig:internal-mu0-sigma0}
\end{figure}

\begin{figure}[!htbp]
	\centering
	\subfloat[]{\includegraphics[scale=0.31]{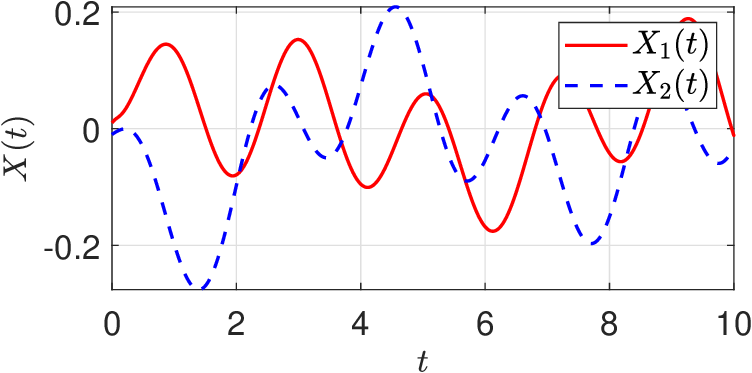}}   \hspace{0.5cm}
	\subfloat[]{\includegraphics[scale=0.31]{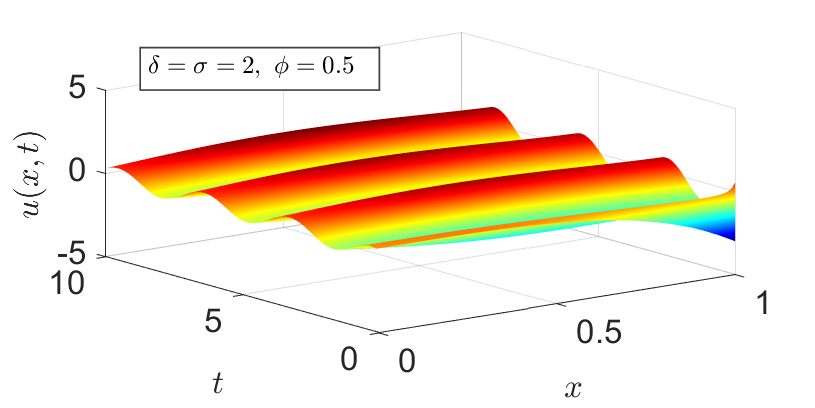}}\\
	\subfloat[]{\includegraphics[scale=0.31]{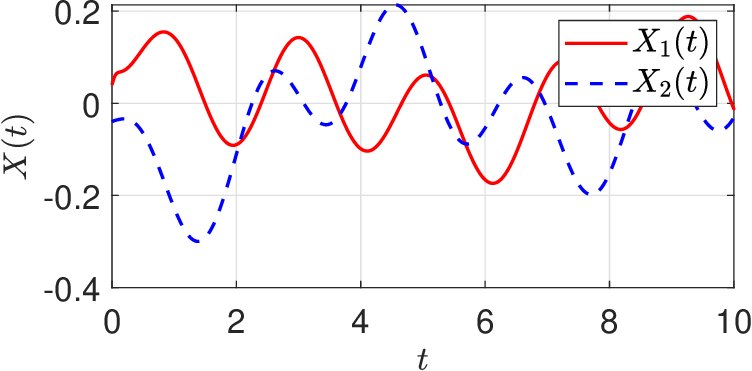}}   \hspace{0.5cm}
	\subfloat[]{\includegraphics[scale=0.31]{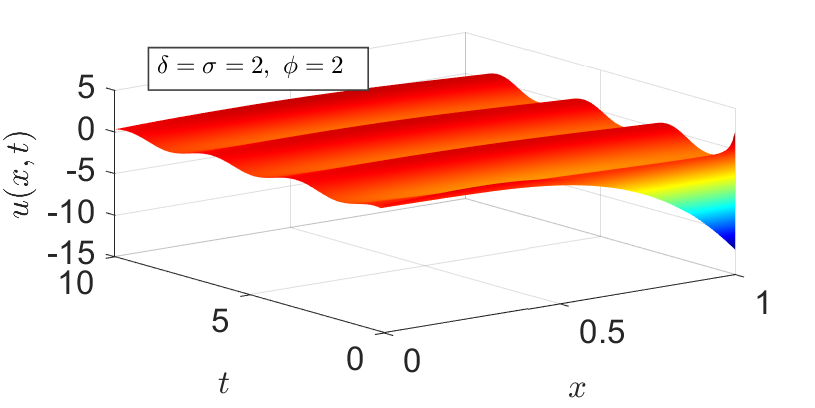}}
	\caption{Evolution of  states    in closed loop with different initial data and the same  internal and in-domain disturbances for $\delta=\sigma=2$: (a) ODE state when $\phi=0.5$; (b) PDE state when $\phi=0.5$; (c) ODE state when $\phi=2$; (d) PDE state when $\phi=2$.}
	\label{fig:internal-mu2-sigma2}
\end{figure}

\begin{figure}[!htbp]
	\centering
	\subfloat[]{\includegraphics[scale=0.31]{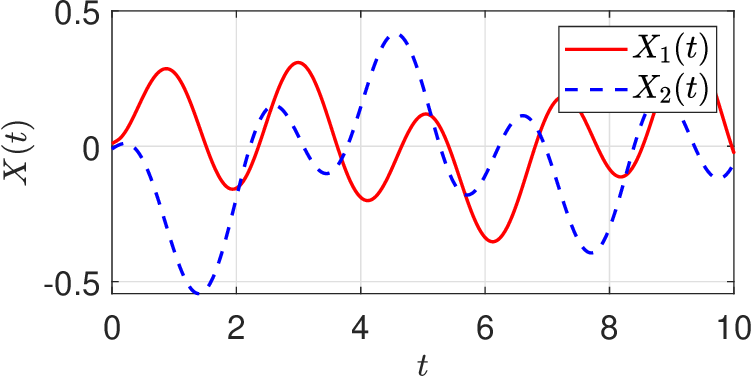}}   \hspace{0.5cm}
	\subfloat[]{\includegraphics[scale=0.31]{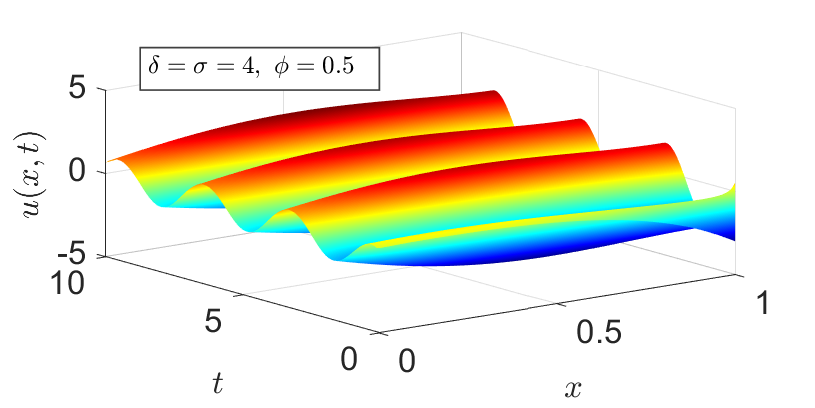}}\\
	\subfloat[]{\includegraphics[scale=0.31]{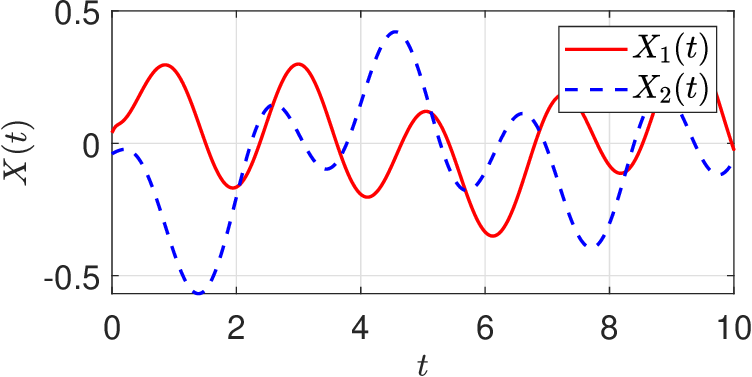}}   \hspace{0.5cm}
	\subfloat[]{\includegraphics[scale=0.31]{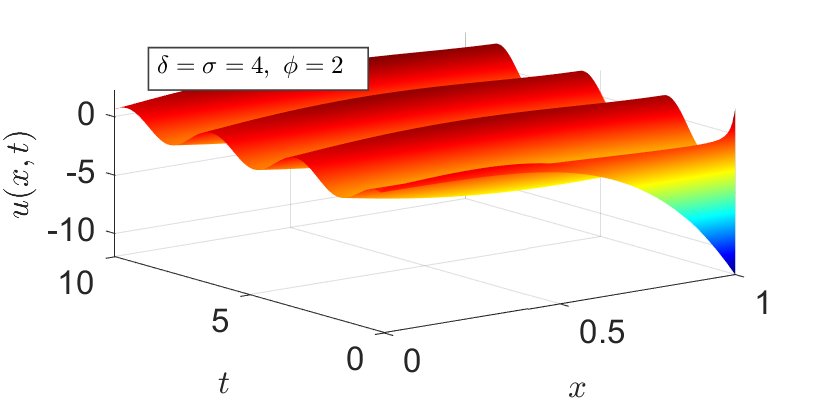}}
	\caption{Evolution of  states    in  closed loop with different initial data and the same internal and in-domain disturbances   for $\delta=\sigma=4$: (a) ODE state when $\phi=0.5$; (b) PDE state when $\phi=0.5$; (c) ODE state when $\phi=2$; (d) PDE state when $\phi=2$.}
	\label{fig:internal-mu4-sigma4}
\end{figure}

\begin{figure}[!htbp]
	\centering
	\subfloat[]{\includegraphics[scale=0.31]{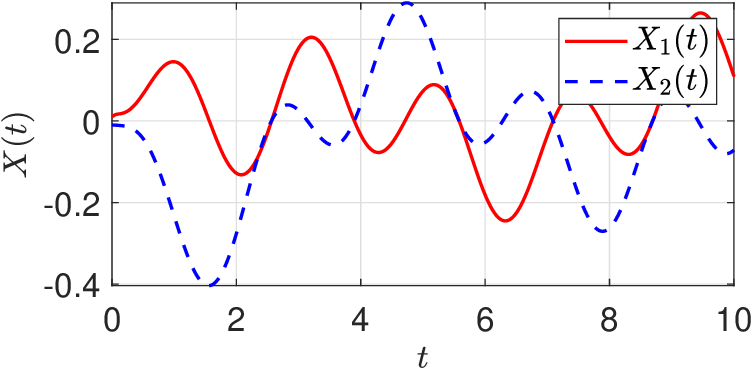}}   \hspace{0.5cm}
	\subfloat[]{\includegraphics[scale=0.31]{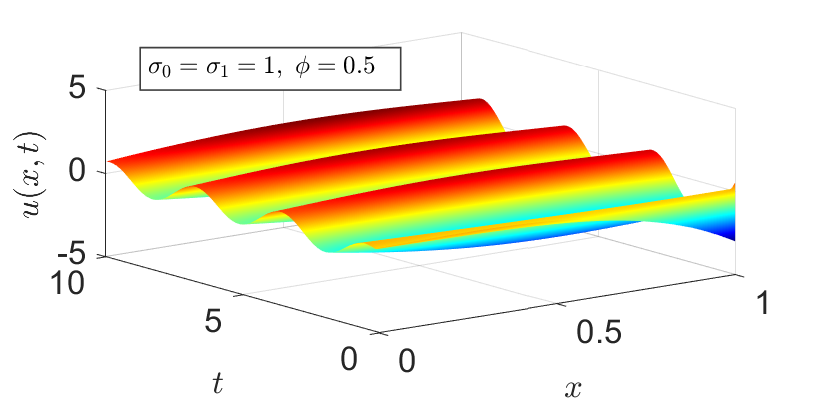}}\\
	\subfloat[]{\includegraphics[scale=0.31]{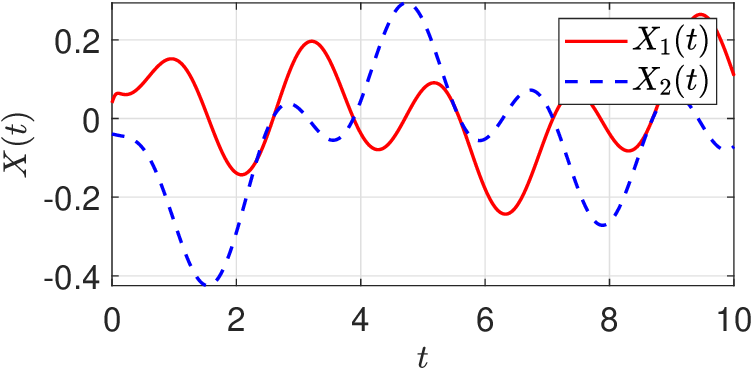}}   \hspace{0.5cm}
	\subfloat[]{\includegraphics[scale=0.31]{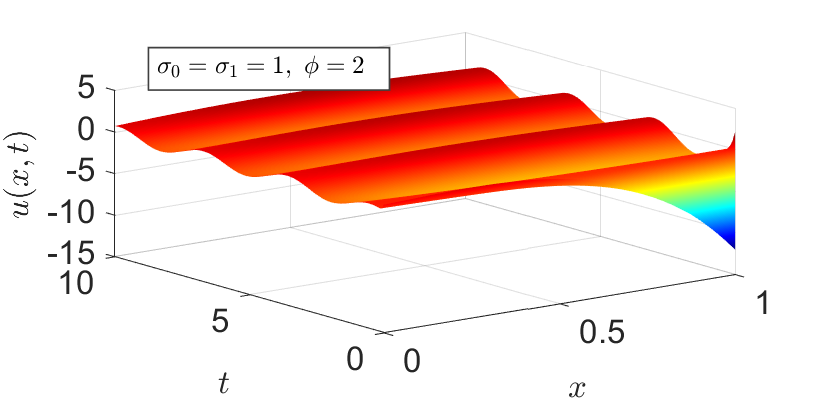}}
	\caption{Evolution of  states    in   closed loop  with different initial data and    the same  boundary disturbances     for  $\sigma_0=\sigma_1=1$: (a) ODE state when $\phi=0.5$; (b) PDE state when $\phi=0.5$; (c) ODE state when $\phi=2$; (d) PDE state when $\phi=2$.}
	\label{fig:boundary-sigma1}
\end{figure}

\begin{figure}[!htbp]
	\centering
	\subfloat[]{\includegraphics[scale=0.31]{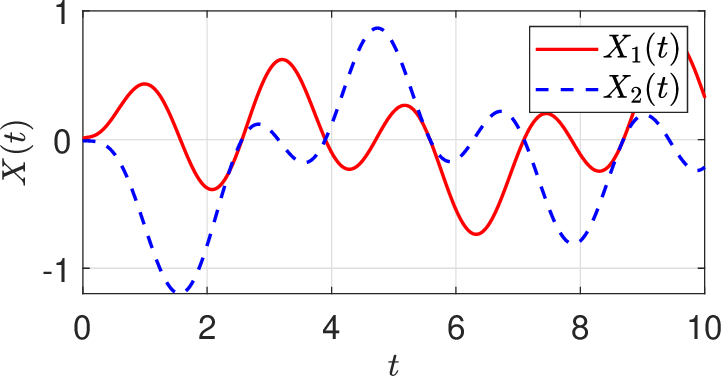}}   \hspace{0.5cm}
	\subfloat[]{\includegraphics[scale=0.31]{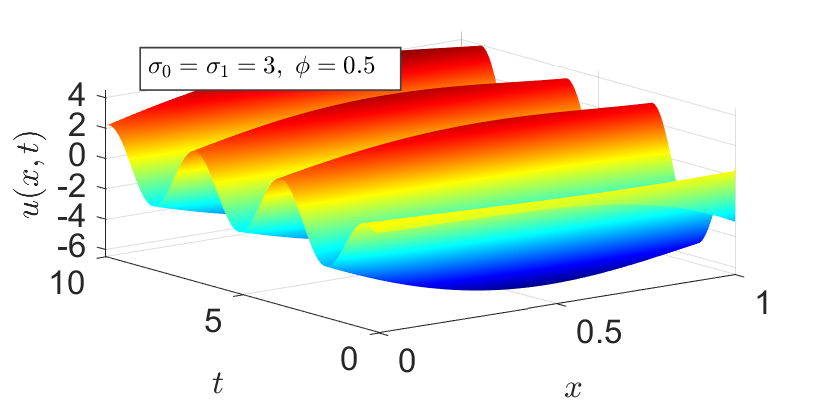}}\\
	\subfloat[]{\includegraphics[scale=0.31]{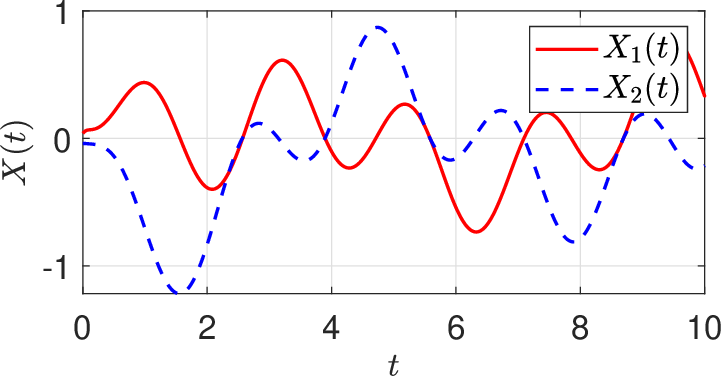}}   \hspace{0.5cm}
	\subfloat[]{\includegraphics[scale=0.31]{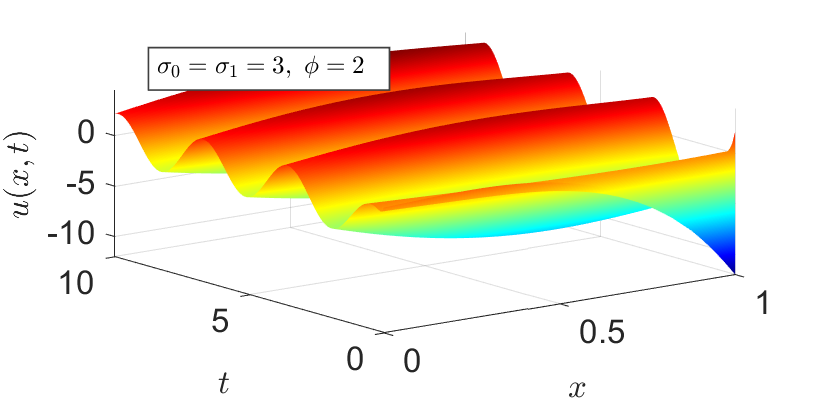}}
	\caption{Evolution of  states    in  closed loop with different initial data and  the same boundary disturbances for $\sigma_0=\sigma_1=3$: (a) ODE state when $\phi=0.5$; (b) PDE state when $\phi=0.5$; (c) ODE state when $\phi=2$; (d) PDE state when $\phi=2$.}
	\label{fig:boundary-sigma3}
\end{figure}

\begin{figure}[!htbp]
	\centering
	\subfloat[]{\includegraphics[scale=0.23]{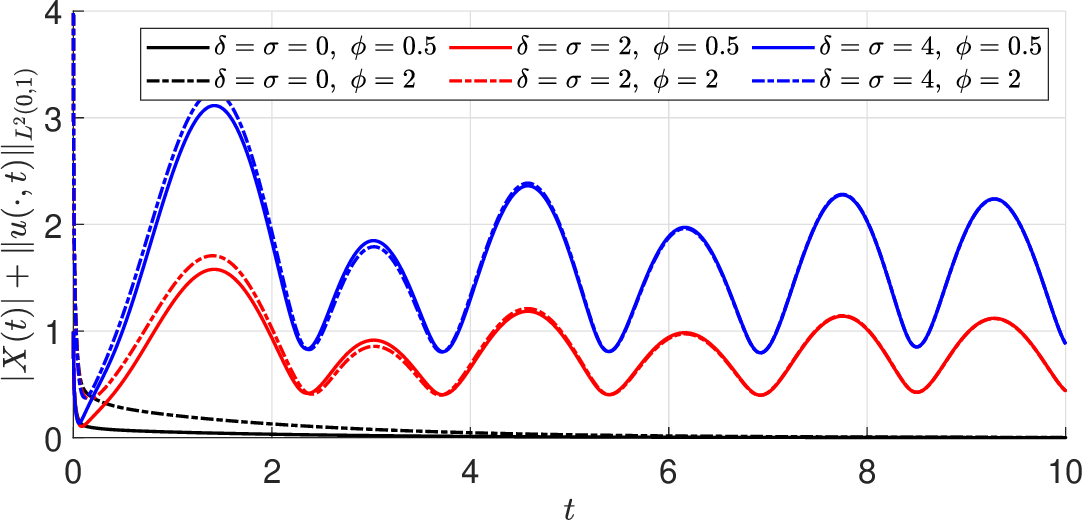}}   \hspace{0.5cm}
	\subfloat[]{\includegraphics[scale=0.23]{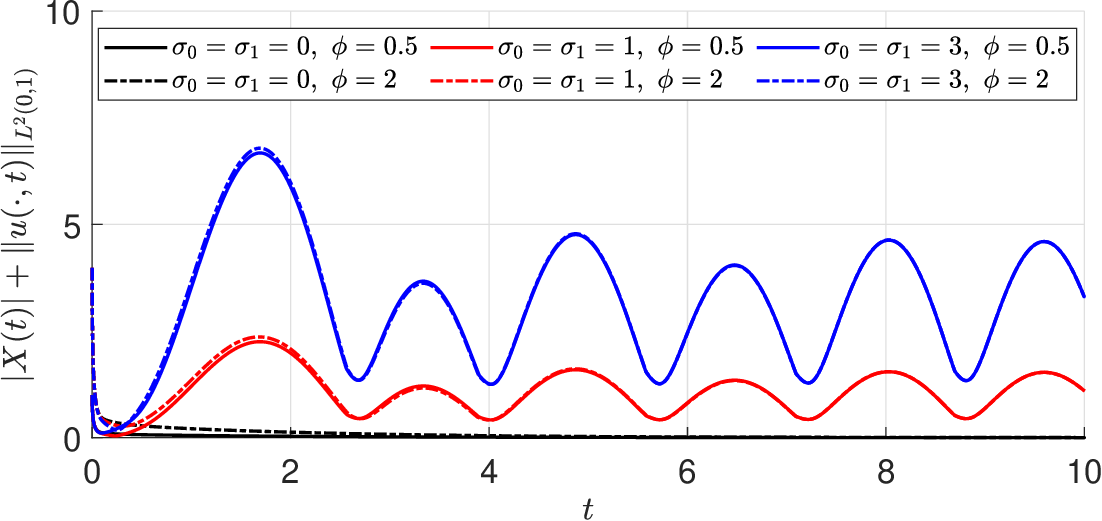}}
	\caption{Evolution of  states' norm   in closed loop with    different    initial data and disturbances: (a) states' norm  when $\phi\in\{0.5,2\}$ and $\delta=\sigma\in\{0,2,4\}$; (b) states' norm   when $\phi\in\{0.5,2\}$ and $\sigma_0=\sigma_1\in\{0,1,3\}$.}
	\label{fig:boundary-norm}
\end{figure}

Figure~\ref{fig:open-loop} shows the evolution of the open-loop system for two sets of initial data, described by $\phi=0.5$ and $\phi=2$. It can be seen that the open-loop system is unstable and that the amplitudes of $X(t)$ and  $u(x,t)$ 
grow rapidly when the magnitude of the initial data becomes larger. In contrast,  Fig.~\ref{fig:internal-mu0-sigma0} shows that, under the composite control law~\eqref{controller-U}, the states of the disturbance-free closed-loop system tend to the equilibrium point. This illustrates the validity of the choice of $l$ and the effectiveness of the proposed composite control law~\eqref{controller-U} in stabilizing the  system~\eqref{original system}.
Figure~\ref{fig:internal-mu2-sigma2} shows that, under the composite control law~\eqref{controller-U}, the states of the closed-loop system with different initial data remain bounded under the same internal and in-domain disturbances with $\delta=\sigma=2$. Especially, for the same  disturbances, the amplitudes of the states become larger when the magnitude of the initial data becomes larger. Figures~\ref{fig:internal-mu2-sigma2} and~\ref{fig:internal-mu4-sigma4} further show that, for the same initial datum, the amplitudes of the states and their ultimate oscillations decrease when the amplitudes of the internal and in-domain disturbances decrease from $\delta=\sigma=4$ to $\delta=\sigma=2$. In accordance with Fig.~\ref{fig:internal-mu2-sigma2} and~\ref{fig:internal-mu4-sigma4}, the evolutions of the corresponding norms are presented in Fig.~\ref{fig:boundary-norm}(a), where the solid and dashed curves correspond to $\phi=0.5$ and $\phi=2$, respectively. At the initial stage, the norm responses are mainly affected by the initial data. As time evolves, the effect of the initial data gradually decreases, and the norm responses are mainly governed by the   disturbances.
In particular, for the same initial datum, the amplitudes of the ultimate norm oscillations decrease when the amplitudes of the internal and in-domain disturbances decrease.
These results depict the ISS property of the closed-loop system w.r.t.   internal and in-domain disturbances $D(t)$ and $f(x,t)$.
Figures~\ref{fig:boundary-sigma1} and~\ref{fig:boundary-sigma3} show that, in the absence of internal and in-domain disturbances, namely, when $\delta=\sigma=0$, the states of the closed-loop system remain bounded under the boundary disturbances with $\sigma_0=\sigma_1=1$ and $\sigma_0=\sigma_1=3$, respectively. Especially, for the same initial datum, the amplitudes of the states and their ultimate oscillations decrease when the amplitudes of the boundary disturbances decrease from $\sigma_0=\sigma_1=3$ to $\sigma_0=\sigma_1=1$.  In accordance with Fig.~\ref{fig:boundary-sigma1} and Fig.~\ref{fig:boundary-sigma3}, the evolutions of the corresponding norms are presented in Fig.~\ref{fig:boundary-norm}(b), where the solid and dashed curves correspond to $\phi=0.5$ and $\phi=2$, respectively. These results depict the ISS property of the closed-loop system w.r.t. boundary disturbances $d_0(t)$ and $d_1(t)$.
Overall, the numerical results in Figs.~\ref{fig:open-loop}-\ref{fig:boundary-norm} show that, with the proposed composite control law~\eqref{controller-U} and  the constant gain $l$,   the stability of the disturbance-free closed-loop system~\eqref{original system} and  the ISS  in the $L^2$-norm  of the closed-loop system with  different external disturbances  can be ensured.


\section{Conclusions}\label{conclusion}
In this paper,   a novel composite controller was designed for the boundary stabilization problem of a coupled ODE-parabolic PDE system with time-varying coefficients in both the ODE and PDE plants, where disturbances enter the ODE plant, the PDE domain, and the Dirichlet-Robin boundary of the PDE.
More specifically, based on an analytic  gain function and a time-varying kernel function,  the composite boundary controller was designed through a composite backstepping transformation, which enables the analytic gain and the time-varying kernel  to be solved successively and avoids directly solving coupled hyperbolic-parabolic kernel equations.
Due to the time-varying coefficients, the ODE-PDE coupling terms, and the  Dirichlet boundary disturbance, the ISS analysis becomes more challenging when the Lyapunov method is applied. To overcome these difficulties,  a square root matrix and a superlinear function were employed to construct  an NQLF   for the ODE and a GLF for the PDE, and the closed-loop stability  was established in the $L^2$-norm.  This provides a new analytical tool for the ISS analysis of coupled ODE-parabolic PDE systems with time-varying coefficient and multiple disturbances, including Dirichlet boundary disturbances.
In the future,  we will consider output-feedback stabilization with only partial measurements of the coupled ODE-PDE system, as well as the extension of the proposed method to nonlinear coupled systems.


\section*{Appendix}

\subsection{Proof of Theorem~\ref{kernel-solution-exitence}}\label{kernel-exitence-proof}
We prove Proposition~\ref{kernel-solution-exitence} by the method of successive approximations; see, e.g.,~\cite{Krstic2008book}. Let
$
\xi:=x+y,  \eta:=x-y,
$
and
$
\Omega_0:=\{(\xi,\eta):0\le \eta\le 1,\ \eta\le \xi\le 2-\eta\}.$
Define
$
G(\xi,\eta,t)
:=
k\left(\frac{\xi+\eta}{2},\frac{\xi-\eta}{2},t\right).
$
By the definition of $\Gamma$ (see~\eqref{Gamma-def}), the kernel equation~\eqref{kernel-k} takes the form over $\Omega_0\times \mathbb{R}_{\geq 0}$: \begin{subequations} \label{G-eq}
	\begin{align}
		G_{\xi\eta}(\xi,\eta,t) & = \frac14 q(t)G(\xi,\eta,t) +\frac14 G_t(\xi,\eta,t), \label{eq:G1}\\ G(\xi,0,t) & = -\frac{\xi}{4}q(t), \label{eq:G2}\\ G_\xi(\xi,\xi,t) & =G_\eta(\xi,\xi,t)+ lG(\xi,\xi,t)+\rho(\xi,t) -\int_0^\xi G(\xi+s,\xi-s,t)\rho(s,t)\,\mathrm ds, \label{eq:G3} 
	\end{align}
\end{subequations}
where  $q(t):=c(t)+\mu(t)$ and $
\rho(x,t):=-\theta(x,t)B(t)$.
Then, after direct calculations, the equation~\eqref{G-eq} is transformed into the following integral equation:
\begin{align}\label{eq:G-integral}
G(\xi,\eta,t)=G_0(\xi,\eta,t)+  \mathcal T[G](\xi,\eta,t),
\end{align}
where
\begin{align*}
G_0(\xi,\eta,t)
:=	
\int_0^\eta e^{l(\tau-\eta)}
\left(
-\frac12 q(t)-\rho(\tau,t)
\right)\,\mathrm d\tau -\frac14q(t)(\xi-\eta)
\end{align*}
and
\begin{align*}
 \mathcal T[G](\xi,\eta,t) 
&:=
\int_0^\eta e^{l(\tau-\eta)}
\Big(
\frac12q(t)\int_0^\tau G(\tau,s,t)\,\mathrm ds  +\int_0^\tau \!\!G(\tau+s,\tau-s,t)\rho(s,t)\,\mathrm ds
+\frac12\int_0^\tau \!\!G_t(\tau,s,t)\,\mathrm ds	
\Big)\,\mathrm d\tau
\nonumber\\
& \quad
+\frac14q(t)\!\int_\eta^\xi\!\!\int_0^\eta
\!G(\tau,s,t)\,\mathrm ds\,\mathrm d\tau +\frac14\int_\eta^\xi\!\!\int_0^\eta
\!G_t(\tau,s,t)\,\mathrm ds\,\mathrm d\tau .
\end{align*}

Define
$	\Phi_0:=G_0$,
and
$	\Phi_{n+1}:=\mathcal T[\Phi_n]$ for
$n\in\mathbb{N}_0$.
Let
$
A_n:=\sum_{j=0}^{n}\Phi_j.
$
By the linearity of $\mathcal T$,
$	A_{n+1}
=
G_0+\mathcal T[A_n]$.
Hence, if the sequence \(\{A_n(\xi,\eta,t)\}\) converges uniformly w.r.t. $(\xi,\eta,t)$ when $n\rightarrow\infty$, then
$G:=\lim_{n\rightarrow\infty}A_n=\sum_{n=0}^{\infty}\Phi_n
$ is a solution to the integral equation
~\eqref{eq:G-integral}.  Therefore, it suffices to establish the uniform convergence of the series $\sum_{n=0}^{\infty}\Phi_n$ w.r.t. $(\xi,\eta,t)$.
Since the integral operator
\(\mathcal T\) contains the time derivative of its argument, the convergence
analysis is carried out for the time derivatives of the iterates. In particular,
we derive estimates on \(\partial_t^r\Phi_n\) for all \(r\in\mathbb N_0\), which also
provide the required regularity of $\Phi_n$ w.r.t. \(t\).

For the initial term, direct differentiation w.r.t. \(t\) gives, for any
\(r\in\mathbb N_0\),
\begin{align*}
	\partial_t^r\Phi_0(\xi,\eta,t) = -\frac14\frac{\mathrm d^r}{\mathrm dt^r}q(t)(\xi-\eta) +\int_0^\eta e^{l(\tau-\eta)}
(
-\frac12\frac{\mathrm d^r}{\mathrm dt^r}q(t)-\partial_t^r\rho(\tau,t)
)\,\mathrm d\tau.
\end{align*}
By the Leibniz's formula,~\eqref{theta-partial},~\eqref{eq:analytic-factorial-assumption-3}, and
\(\binom{r}{j}(r-j)!j!=r!\), we have
\begin{align}	\label{eq:q-rho-majorant}
\left|\partial_t^r\rho(x,t)\right|
&\le
\sum_{j=0}^r
\binom{r}{j}
\left(2ae^b b^{r-j}(r-j)!\right)
\left(\alpha^{j+1}j!\right) 
 \le  \gamma^{2r+2}r!,
\ \forall r\in\mathbb N_0,
\end{align}
where \(\gamma:=\max\{2ae^b,\alpha,b,2\}\). 
By~\eqref{eq:analytic-factorial-assumption-3} and~\eqref{mu-ana}, we have
\begin{align}\label{eq-75}
\sup_{t\in\mathbb{R}_{\geq 0}}
\left|\frac{\mathrm d^r}{\mathrm dt^r} q(t)\right|
&\le
\gamma^{r+2}r!,\ \forall r\in\mathbb{N}_0.
\end{align}
Using~\eqref{eq:q-rho-majorant},~\eqref{eq-75}, and
\(|e^{l(\tau-\eta)}|\le 1\) for \(0\le\tau\le\eta\), 
we obtain
\begin{align}\label{n=0}
\!\!\! \sup_{t\in\mathbb{R}_{\geq 0}}
\left|\partial_t^r\Phi_0(\xi,\eta,t)\right|
\!  \le
\gamma^{2r+3}(r+1)!(\xi+\eta),
\ \forall r\in\mathbb{N}_0.\!
\end{align}
In particular, for \(r=0\), it holds that
$
|\Phi_0(\xi,\eta,t)|\le \gamma^3(\xi+\eta).
$

We claim that, for all \(n\in\mathbb{N}_0\),
\begin{align}
\sup_{t\in\mathbb{R}_{\geq 0}}
\left|\partial_t^r \Phi_n(\xi,\eta,t)\right|\le
\gamma^{2n+2r+3}
\frac{(n+r+1)!}{(n!)^2}
(\xi\eta)^n(\xi+\eta),\ \forall r\in\mathbb{N}_0.
\label{eq:claim-estimate}
\end{align}
For \(n=0\),~\eqref{eq:claim-estimate} reduces to~\eqref{n=0}. Assume that~\eqref{eq:claim-estimate} holds for  any fixed \(n\in\mathbb N_0\). We prove that it
also holds for \(n+1\). Indeed, for any \(r\in\mathbb N_0\), differentiating
\(\Phi_{n+1}\) w.r.t. \(t\) gives
\begin{align}\label{phi_n+1}
\sup_{t \in\mathbb{R}_{\geq 0}} |\partial_t^r\Phi_{n+1}(\xi,\eta,t)|
&\le  \frac12
\int_0^\eta 	\int_0^\tau
\sup_{t \in\mathbb{R}_{\geq 0}} |\partial_t^r (q(t)\Phi_n(\tau,s,t))|\,\mathrm ds
\,\mathrm d\tau +
\frac12
\int_0^\eta
\int_0^\tau
\sup_{t \in\mathbb{R}_{\geq 0}} |\partial_t^{r+1}\Phi_n(\tau,s,t)|\,\mathrm ds
\,\mathrm d\tau \notag\\
&\quad+
\int_0^\eta 	\int_0^\tau
\sup_{t \in\mathbb{R}_{\geq 0}} |\partial_t^{r}(\Phi_n(\tau+s,\tau-s,t)\rho(s,t))|\,\mathrm ds
\,\mathrm d\tau \notag\\
&\quad+
\frac14
\int_\eta^\xi\int_0^\eta
\sup_{t \in\mathbb{R}_{\geq 0}}
|\partial_t^r (q(t)\Phi_n(\tau,s,t))|
\,\mathrm ds\,\mathrm d\tau  +
\frac14
\int_\eta^\xi\int_0^\eta
\sup_{t \in\mathbb{R}_{\geq 0}} |\partial_t^{r+1}\Phi_n(\tau,s,t)|
\,\mathrm ds\,\mathrm d\tau.
\end{align}
Using
the Leibniz's formula,~\eqref{eq:claim-estimate}, and~\eqref{eq:q-rho-majorant}, and then applying   
$	\sum_{j=0}^r\binom{r}{j} j!(n+r+1-j)!
\le \frac{(n+r+2)!}{n+1}$,
we obtain
\begin{align}	\label{eq:q-term-estimate}
 \sup_{t \in\mathbb{R}_{\geq 0}}
\left|
\partial_t^r
\left(q(t)\Phi_n(\tau,s,t)\right)
\right| 
 \le
\gamma^{2n+2r+5}
\frac{(n+r+2)!}{(n+1)(n!)^2}
(\tau s)^n(\tau+s).
\end{align}
Moreover, by~\eqref{eq:claim-estimate} with \(r\) replaced by \(r+1\), we have
\begin{align}
\sup_{t \in\mathbb{R}_{\geq 0}}\!
\left|
\partial_t^{r+1}\Phi_n(\tau,s,t)
\right| \le
\gamma^{2n+2r+5}
\frac{(n+r+2)!}{(n!)^2}
(\tau s)^n(\tau+s).
\label{eq:t-term-estimate}
\end{align}
Using the Leibniz's formula,~\eqref{eq:claim-estimate}, and~\eqref{eq:q-rho-majorant}, we obtain
\begin{align}
 \sup_{t \in\mathbb{R}_{\geq 0}}
\left|
\partial_t^r
\left(
\Phi_n(\tau+s,\tau-s,t)\rho(s,t)
\right)
\right| \le
\gamma^{2n+2r+5}
\frac{(n+r+2)!}{(n+1)(n!)^2}
\big((\tau+s)(\tau-s)\big)^n(2\tau).
\label{eq:rho-term-estimate}
\end{align}
Putting~\eqref{eq:q-term-estimate},~\eqref{eq:t-term-estimate}, and~\eqref{eq:rho-term-estimate} into~\eqref{phi_n+1}, we deduce that
\begin{align}\label{phi_n+12}
 \sup_{t \in\mathbb{R}_{\geq 0}} |\partial_t^r\Phi_{n+1}(\xi,\eta,t)| 
 	&\le  \frac12 	 \gamma^{2n+2r+5}
\frac{(n+r+2)!}{(n+1)(n!)^2}
\int_0^\eta 	\int_0^\tau
(\tau s)^n(\tau+s)\,\mathrm ds
\,\mathrm d\tau \notag\\
\!\!\!	&\quad+
\frac12	 \gamma^{2n+2r+5}
\frac{(n+r+2)!}{(n!)^2}
\int_0^\eta
\int_0^\tau
(\tau s)^n(\tau+s)\,\mathrm ds
\,\mathrm d\tau \notag\\
\!\!\!	&\quad+  \gamma^{2n+2r+5}
\frac{(n+r+2)!}{(n+1)(n!)^2}
\int_0^\eta 	\int_0^\tau
\tau^{2n}(2\tau)\,\mathrm ds
\,\mathrm d\tau \notag\\
\!\!\!	&\quad+
\frac14   \gamma^{2n+2r+5}
\frac{(n+r+2)!}{(n+1)(n!)^2}
\int_\eta^\xi\int_0^\eta
(\tau s)^n(\tau+s)
\,\mathrm ds\,\mathrm d\tau \notag\\
\!\!\!	&\quad+
\frac14	 \gamma^{2n+2r+5}
\frac{(n+r+2)!}{(n!)^2}
\int_\eta^\xi\int_0^\eta
(\tau s)^n(\tau+s)
\,\mathrm ds\,\mathrm d\tau.
\end{align}
For \((\xi,\eta)\in\Omega_0\), the following integral estimates hold:
\begin{subequations}
\begin{align}
\int_0^\eta 	\int_0^\tau
(\tau s)^n(\tau+s)\,\mathrm ds
\,\mathrm d\tau
&\le
\frac{(\xi\eta)^{n+1}(\xi+\eta)}{2(n+1)(n+2)},\label{eq:An}\\
\int_0^\eta\int_0^\tau 2\tau^{2n+1}\,\mathrm ds\,\mathrm d\tau
&\le
\frac{(\xi\eta)^{n+1}(\xi+\eta)}{2n+3},
\label{eq:Bn}\\
\int_\eta^\xi\int_0^\eta(\tau s)^n(\tau+s)\,\mathrm ds\,\mathrm d\tau
&\le
\frac{(\xi\eta)^{n+1}(\xi+\eta)}
{(n+1)(n+2)}.
\label{eq:Cn}
\end{align}
\end{subequations}
Substituting~\eqref{eq:An},~\eqref{eq:Bn}, and~\eqref{eq:Cn} into~\eqref{phi_n+12} yields
$\sup_{t\geq0}|\partial_t^r\Phi_{n+1}(\xi,\eta,t)|
\! \le \!
\gamma^{2n+2r+5}
\frac{(n+r+2)!}{(n+1)!}
\frac{(\xi\eta)^{n+1}}{(n+1)!}(\xi+\eta)$  
for any $ r\in\mathbb{N}_0$.
Thus, by induction,~\eqref{eq:claim-estimate} holds for $n+1$, and hence for all
\(n\in\mathbb{N}_0\).

Taking \(r=0\) in~\eqref{eq:claim-estimate}, we get
$	|\Phi_n|
\le
\gamma^3(\xi+\eta)(n+1)\frac{(\gamma^2\xi\eta)^n}{n!}$.
Moreover, since \((\xi,\eta)\in\Omega_0\), we have \(0\le\xi+\eta\le 2\) and
\(\xi\eta\le1\). Hence
$	 |G|
\le \sum_{n=0}^{\infty}|\Phi_n|
\le
2\gamma^3(1+\gamma^2)e^{\gamma^2}$.
Thus, the series \(\sum_{n=0}^{\infty}\Phi_n\) converges absolutely
and uniformly w.r.t. $(\xi,\eta,t)$ on \(\Omega_0\times\mathbb{R}_{\geq 0}\), and 	\(G=\sum_{n=0}^{\infty}\Phi_n\) is a solution to~\eqref{eq:G-integral}. Moreover, $G$ is continuous in $(\xi,\eta)$. Since $q$ and $\rho$ are analytical, it follows from~\eqref{eq:G-integral} that \(G\) is of class \(C^2\)   w.r.t. $(\xi,\eta)$. Furthermore, for each \(r\ge1\), the estimate~\eqref{eq:claim-estimate} also imply the uniform convergence of
\(\sum_{n=0}^{\infty}\partial_t^r\Phi_n\) w.r.t. $(\xi,\eta,t)$. Therefore, by using the change of variables, we deduce that the kernel~\eqref{kernel-k} admits a solution $k\in C^{\infty}(\mathbb{R}_{\geq 0};C^2(\mathcal{D}_0))$.

In addition, since \(\xi+\eta=2x\) and \(\xi\eta=x^2-y^2\),
using
\(k(x,y,t)=\sum_{n=0}^{\infty}\Phi_n(x+y,x-y,t)\),
the estimate~\eqref{eq:claim-estimate}  yields
\begin{align}\label{eq-109}
\!\! \sup_{t\in\mathbb{R}_{\ge0}}\!\left|\partial_t^r k(x,y,t)\right|    \le 
2x\gamma^{2r+3}\!\sum_{n=0}^{\infty} 
\frac{(n+r+1)!}{(n!)^2}\!
\left(\gamma^2(x^2-y^2)\right)^{\!n} =\!   2x\gamma^{2r+3}\! \frac{\mathrm d^{r\! +\! 1}}{\mathrm dz^{r\! +\! 1}}(z^{r\! +\! 1}\! e^z)\Big|_{z=\gamma^2(x^2\! -\! y^2)}\! .\!
\end{align}
Note that the Leibniz's formula implies that
\begin{align}\label{dhe}
\frac{\mathrm d^{r+1}}{\mathrm dz^{r+1}}
\left(z^{r+1}e^z\right) 
\le
(r+1)!e^z(1+z)^{r+1},\forall z\in \mathbb{R}.
\end{align}
Putting~\eqref{dhe} into~\eqref{eq-109}, 
we get a uniform upper bound for $	\sup_{t\in\mathbb R_{\ge0}}
\left|\partial_t^i k(x,y,t)\right|$ w.r.t. $t$ in Theorem~\ref{kernel-solution-exitence}.

It remains to prove the  uniqueness of solutions $G$ to the integral equation~\eqref{eq:G-integral}. Suppose that \(G^{(1)}\) and \(G^{(2)}\) are two
solutions to~\eqref{eq:G-integral}. Define
$
\widetilde G:=G^{(1)}-G^{(2)}.
$
Then \(\widetilde G=\mathcal T[\widetilde G]\). Define $\Psi_0:=\widetilde G$ and $\Psi_{n+1}:=\mathcal T[\Psi_n]$ for $n\in\mathbb{N}_0$. We claim that
$
\Psi_n=\widetilde G
$ for all $ n\in\mathbb N_0$.
Indeed, this is true for \(n=0\). Moreover, if \(\Psi_n=\widetilde G\), then
$
\Psi_{n+1}=\mathcal T[\Psi_n]=\mathcal T[\widetilde G]=\widetilde G,
$
where \(\widetilde G=\mathcal T[\widetilde G]\) has been used. Hence the claim follows by induction.
Using~\eqref{eq:claim-estimate} and the definition of $G$,
\(\Psi_0 \) satisfies
\begin{align}\label{partial-H}
\!\!\!\!	\sup_{t\in\mathbb{R}_{\geq0}}
\left|\partial_t^r \Psi_0(\xi,\eta,t)\right|
&\leq \!2\gamma^{2r+3}(\xi+\eta)\!\frac{\mathrm d^{r+1}}{\mathrm d{z}^{r+1}}
\!\left(z^{r+1}e^{z}\right)\Big|_{z=\gamma^2\xi\eta}.\!
\end{align}
Note that 
$
(1+\gamma^2\xi\eta)^{r+1}\le (1+\gamma^2)^{r+1}.
$
Applying~\eqref{dhe} to~\eqref{partial-H}  {and letting 	\(\hat\gamma\ge\gamma(1+\gamma^2)\)}, we have
 \begin{align*}
 	\sup_{t\in\mathbb{R}_{\geq0}}
\left|\partial_t^r \Psi_0(\xi,\eta,t)\right|   \le
2\gamma^{2r+3}(r+1)!
(1 + \gamma^2)^{r+1}
(\xi + !\eta)e^{\gamma^2\xi\eta}
\leq 2\hat\gamma^{2r+3}(r+1)!(\xi+\eta)e^{\hat\gamma^2\xi\eta}.
 \end{align*}

Analogous to~\eqref{eq:claim-estimate}, we deduce that, for every \(n\in\mathbb N_0\),
\begin{align}\label{unique-G}
 \sup_{t\in\mathbb{R}_{\geq 0}} |\partial_t^r \Psi_n(\xi,\eta,t)|  \le 
2e^{\hat{\gamma}^2\xi\eta}
\hat{\gamma}^{2n+2r+3}
\frac{(n\!+r\!+1)!}{(n!)^2}
(\xi\eta)^n(\xi+\eta),\ \! \forall r\in\mathbb N_0.
\end{align}

Taking \(r=0\) in~\eqref{unique-G}, we obtain
$	|\Psi_{n}|
\le
2e^{\hat{\gamma}^2\xi\eta}
\hat{\gamma}^{2n+3}
(n+1)\frac{(\xi\eta)^n}{n!}
(\xi+\eta)$ for all
$n\in\mathbb N_0$.
Since \(\Psi_n=\widetilde G\) for all \(n\in\mathbb{N}_0\) and the left-hand side
above tends to zero as \(n\to\infty\), it follows that
\(\widetilde G=0\). Hence
\(G^{(1)}=G^{(2)}\), which proves the uniqueness of the solution $G$ to~\eqref{eq:G-integral}.
\subsection{Proof of Proposition~\ref{target-sytem-prop}}\label{proof-of-target}
We derive the target system~\eqref{target-system} from the composite transformation~\eqref{original system} in two steps.

\textbf{Step 1:} Transform $(X,u)$ to $(X,v)$ via~\eqref{transfor-1}. Specially, we show that system~\eqref{original system} can be transformed into\begin{subequations}\label{v-system}
\begin{align}
\!\!\!\dot X(t)&= F(t) X(t)+B(t)v(0,t) +D(t), \ t\in \mathbb{R}_{>0},\label{v-ODE}\\
\!\!\!v_t(x,t)&=v_{xx}(x,t)+c(t)v(x,t)-\theta(x,t)B(t)v(0,t) +f(x,t)-\theta(x,t)D(t),\ (x,t)\in Q_\infty,\label{v-eq1}\\
\!\!\!v_x(0,t)	&=l v(0,t)+d_0(t),\ t\in \mathbb{R}_{>0},\label{v-eq2}\\
\!\!\!v(1,t)&=U(t)+d_1(t)-\theta(1,t)X(t),\ t\in \mathbb{R}_{>0},\label{v-eq3}\\
\!\!\!X(0)&=X_0,\ v(x,0)=v_0(x), \  x\in(0,1),\label{initial-vale}
\end{align}
\end{subequations}
where $v_0(x):=u_0(x)-\theta(x,0)X_0$.

By the transformation~\eqref{transfor-1} and~\eqref{theta-boundary1},  we have
$	\dot X(t)
=\big(A(t)+B(t)K(t)\big)X(t)+B(t)v(0,t)+D(t)$,
which implies that~\eqref{v-ODE} holds.
Differentiating~\eqref{transfor-1} twice w.r.t. \(x\) and once w.r.t. \(t\), and then using~\eqref{original system} and~\eqref{transfor-1}, we obtain
\begin{align*}
 	v_t(x,t)-v_{xx}(x,t)-c(t)v(x,t) 
&=-\theta(x,t)B(t)v(0,t)+f(x,t)-\theta(x,t)D(t)\\
&\quad +\Big(\theta_{xx}(x,t)-\theta_t(x,t)
-\theta(x,t)(A(t)+B(t)\theta(0,t)) +c(t)\theta(x,t)\Big)X(t)\notag\\
&=-\theta(x,t)B(t)v(0,t)+f(x,t)-\theta(x,t)D(t),
\end{align*}
which implies that~\eqref{v-eq1} holds.
For the boundary condition at \(x=0\), differentiating~\eqref{transfor-1} once w.r.t. \(x\) and using~\eqref{PDE2},~\eqref{transfor-1}, and~\eqref{theta-boundary2}, we have
$	v_x(0,t)
= l v(0,t)+\big(l\theta(0,t)+C(t)-\theta_x(0,t)\big)X(t)+d_0(t)
= l v(0,t)+d_0(t)$.
At \(x=1\), using~\eqref{transfor-1} and~\eqref{PDE3}, we have
$
v(1,t)=U(t)+d_1(t)-\theta(1,t)X(t).
$
Then we get the boundary conditions of $v$, i.e., ~\eqref{v-eq2} and~\eqref{v-eq3}.

Moreover, combining~\eqref{PDE4} with~\eqref{transfor-1}, we obtain the initial condition~\eqref{initial-vale}. Therefore, under the transformation~\eqref{transfor-1}, the   state $(X,u)$ is mapped into $(X,v)$, which satisfies~\eqref{v-system}.

\textbf{Step 2:} Transform  $(X,v)$ to $(X,w)$  via~\eqref{transfor-2}. Substituting~\eqref{transfor-2} into~\eqref{v-ODE}, we directly obtain~\eqref{w-ODE}. Then, proceeding as in Step~1, we differentiate~\eqref{transfor-2} w.r.t. \(x\) and \(t\),
apply integration by parts twice, and use~\eqref{transfor-2} again to obtain
\begin{align}\label{eq-126}
\!\!\! 	w_t(x,t)-w_{xx}(x,t)+\mu(t)w(x,t) 
&=\Big(c(t)+\mu(t)+2\frac{\mathrm d}{\mathrm dx}(k(x,x,t))\Big)v(x,t)\notag\\
\!\!\!&\quad +\int_0^x
\Big(k_{xx}(x,y,t)-k_{yy}(x,y,t)-k_t(x,y,t) -(c(t)+\mu(t))k(x,y,t)\Big)v(y,t)\,\mathrm dy\notag\\
\!\!\!&\quad +\Big(\rho(x,t)+lk(x,0,t)-k_y(x,0,t)  -\int_0^x k(x,y,t)\rho(y,t)\,\mathrm dy\Big)v(0,t)+\psi(x,t)\notag\\
\!\!\!&\quad-\int_0^x k(x,y,t)\psi(y,t)\,\mathrm dy
+k(x,0,t)d_0(t),
\end{align}
where $\rho(x,t):=-\theta(x,t)B(t)$ and
$\psi(x,t):=f(x,t)-\theta(x,t)D(t)$. Then,~\eqref{eq-126} implies that~\eqref{target-2} holds.
For the boundary conditions of $w$, by the boundary conditions of $v$, i.e.,~\eqref{v-eq2} and~\eqref{v-eq3},   the transformation~\eqref{transfor-2}, the composite control law~\eqref{controller-U}, and the kernel function equation~\eqref{kernel-k}, we get
$
w_x(0,t)
=l w(0,t)+d_0(t)
$ and 
$w(1,t)
=U(t)+d_1(t)-\theta(1,t)X(t) -\int_0^1 k(1,y,t)v(y,t)\,\mathrm dy
%
=d_1(t)$.
Then, we get the boundary conditions of $w$, i.e.,~\eqref{target-3} and~\eqref{target-4}.
Moreover, combining~\eqref{initial-vale} with~\eqref{transfor-2}, we obtain the initial condition~\eqref{target-5}. Therefore, under the transformation~\eqref{transfor-2}, the   state $(X,v )$ is mapped into  $(X,w )$, which satisfies~\eqref{target-system}.

\subsection{Proof of Proposition~\ref{equivalent}}\label{proof-of-equivent}
We   prove Proposition~\ref{equivalent} in three steps.

\textbf{Step 1:} Prove  the boundedness of the transformation
$\mathcal{K}_1$ and its inverse. By~\eqref{theta-partial} and~\eqref{transfor-1}, we have
\begin{align}\label{vleqx}
\!\!\!	\|(\mathcal{K}_1[X,u])(\cdot,t)\|_\mathcal{H}
\! 	 \leq\!
(1\!+\!\overline\theta)	\|(X(t),\!u(\cdot,\!t))\|_{\mathcal{H}},\!\ t\in\mathbb{R}_{>0},\!
\end{align}
where
$  \overline\theta=2ae^b$. Moreover, the   transformation~\eqref{transfor-1} is explicitly
invertible since
$	u(x,t)=v(x,t)+\theta(x,t)X(t)$,
and hence
\begin{align}\label{u=v+x}
\|(\mathcal{K}^{-1}_1[X,v])(\cdot,t)\|_\mathcal{H}
&\leq
(1+\overline\theta) \|(X(t),v(\cdot,t))\|_{\mathcal{H}}.
\end{align}
Therefore,  the transformation
\(\mathcal{K}_1:(X,u)\mapsto(X,v)\) is invertible and its
inverse is bounded on $\mathcal{H}$.

\textbf{Step 2:} Prove  the boundedness of the transformation
$\mathcal{K}_2$ and its inverse. 
For any fixed $t\in \mathbb{R}_{\geq 0}$, define the Volterra integral operator
\(\mathcal S:L^2((0,1);\mathbb{R})\to L^2((0,1);\mathbb{R})\) by
$	(\mathcal S\varphi)(x)
:=
\int_0^x k(x,y,t)\varphi(y)\,\mathrm dy$ for any $
\varphi\in L^2((0,1);\mathbb{R})$,
and hence
$
(\mathcal K_2[X,v])(\cdot,t):=(X(t),(I-\mathcal S)v(\cdot,t)),
$
where \(I\) denotes the identity operator on \(L^2((0,1);\mathbb{R})\).

First, we  show that $\mathcal{K}_2$ is bounded.
Indeed, for
\(\varphi\in L^2((0,1);\mathbb{R})\), by the Cauchy-Schwarz inequality, we have
\begin{align*}
 \|\mathcal S\varphi\|_{L^2((0,1);\mathbb{R})}^2  =
\int_0^1
\left|
\int_0^x k(x,y,t)\varphi(y)\,\mathrm dy
\right|^2
\mathrm dx
\leq
\frac{\overline k^2}{2}
\|\varphi\|_{L^2((0,1);\mathbb{R})}^2,
\end{align*}
where
$
\overline k
:=
\sup_{(x,y,t)\in\mathcal D}|k(x,y,t)|.
$
Hence \(\mathcal S\varphi\in L^2((0,1);\mathbb{R})\). The linearity of
\(\mathcal S\), and hence that of \(\mathcal K_2\), follows
from the linearity of the integral.  
Consequently,
\begin{align}\label{K2}
\|(\mathcal K_2[X,v])(\cdot,t)\|_{\mathcal{H}} 
&\leq
\left(1+\frac{\overline k}{\sqrt2}\right)
\|(X(t),v(\cdot,t))\|_{\mathcal{H}}.
\end{align}

Next, we prove the invertibility of \(\mathcal K_2\). Note that   \(\mathcal K_{2}\) leaves the $X$
component unchanged, while its $v$ component is given by
\(I-\mathcal S\).   Hence, it suffices to establish the
invertibility of \(I-\mathcal S\). Specifically, we show that
\begin{equation}
(I-\mathcal S)^{-1}
=
I+\sum_{m=1}^{\infty}\mathcal S^m .
\label{K_-1-L2}
\end{equation}
Let \(k_m\) denote the kernel of \(\mathcal S^m\), where \(k_1(x,y,t):=k(x,y,t)\) and
$k_{m+1}(x,y,t)
:=
\int_y^x k(x,z,t)k_m(z,y,t)\,\mathrm dz$ for \(m\geq1\).
We  claim that, for every \((x,y,t)\in\mathcal D\),
\begin{equation}
|k_m(x,y,t)|
\leq
\overline k^m
\frac{(x-y)^{m-1}}{(m-1)!},
\ m\geq1.
\label{eq:iterated-kernel-bound-L2}
\end{equation}
Indeed, the case \(m=1\) follows from the definition of \(\overline k\), since
\(|k_1(x,y,t)|=|k(x,y,t)|\leq \overline k\). If~\eqref{eq:iterated-kernel-bound-L2}
holds for  \(m\geq1\), then
$	|k_{m+1}(x,y,t)|
\!\leq\!\!
\int_y^x\! \!|k(x,z,t)|\!\,|k_m(z,y,t)|\!\,\mathrm dz  \!\leq
\overline k^{m+1}\!
\frac{(x\!-\!y)^m}{m!}$,
proving~\eqref{eq:iterated-kernel-bound-L2} by induction.
It follows from~\eqref{eq:iterated-kernel-bound-L2} and the Young's inequality that, for any $\varphi \in L^2((0,1);\mathbb{R})$,
\begin{align*}
	\|\mathcal S^m\varphi\|_{L^2(0,1)}
\leq
\left(
\int_0^1
\overline k^m
\frac{s^{m-1}}{(m-1)!}\,\mathrm ds
\right)
\|\varphi\|_{L^2((0,1);\mathbb{R})}
=
\frac{\overline k^m}{m!}
\|\varphi\|_{L^2((0,1);\mathbb{R})},
\end{align*}
and hence
\begin{equation}
\|\mathcal S^m\|_{\mathscr L(L^2(0,1))}
\leq
\frac{\overline k^m}{m!},
\ m\geq1.
\label{eq:A-power-bound-L2}
\end{equation}
Therefore, \(\sum_{m=0}^{\infty}\mathcal S^m\) converges in
\(\mathscr L(L^2((0,1);\mathbb{R}))\).
Moreover, for every \(N\geq0\), we have
$
(I-\mathcal S)\sum_{m=0}^{N}\mathcal S^m
=
\left(\sum_{m=0}^{N}\mathcal S^m\right)$ $(I-\mathcal S)
=
I-\mathcal S^{N+1}.
$
By~\eqref{eq:A-power-bound-L2}, letting $N\to\infty$ yields
$(I-\mathcal S)\sum_{m=0}^{\infty}\mathcal S^m
=
\left(\sum_{m=0}^{\infty}\mathcal S^m\right)(I-\mathcal S)
=
I$,
which proves~\eqref{K_-1-L2} and the invertibility of $I-\mathcal S$, and thus of $\mathcal{K}_2$.

Finally, we prove the boundedness of the inverse transformation $\mathcal{K}_2^{-1}$. It follows from~\eqref{K_-1-L2} and~\eqref{eq:A-power-bound-L2} that
\begin{align*}
\|(I-\mathcal S)^{-1}\|_{\mathscr L(L^2((0,1);\mathbb{R}))}
\leq
\sum_{m=0}^{\infty}
\frac{\overline k^m}{m!}
=
e^{\overline k}.
\end{align*}
Then, applying \((I-\mathcal S)^{-1}\) to~\eqref{transfor-2}, we obtain  
\begin{align}
\!\!\!\!\!\|(\mathcal{K}_2^{-1}[X,w])(\cdot,t)\|_{\mathcal{H}}
\! \leq
\left(e^{\overline k}+1\right)
\|(X(t),w(\cdot,t))\|_{\mathcal H}.
\label{k2inverse}
\end{align}

\textbf{Step 3:} Prove  the boundedness of \(\mathcal G \) and \(\mathcal G^{-1}\).  From~\eqref{vleqx} and~\eqref{K2}, we obtain~\eqref{bound-1}.
Moreover, combining~\eqref{u=v+x} and~\eqref{k2inverse} gives~\eqref{bound-2}.
%
Therefore,  \(\mathcal G\) is a bounded
bijective linear operator on \(\mathcal H\), and its inverse
\(\mathcal G^{-1}\) is also bounded.

\end{document}